\documentclass[12pt,a4paper,final]{iopart}

\usepackage{iopams}  
\usepackage{graphicx}
\usepackage[breaklinks=true,colorlinks=true,linkcolor=blue,urlcolor=blue,citecolor=blue]{hyperref}
\usepackage{amssymb}
\newcommand{\R}{{\mathbb R}}

\newcommand{\N}{{\mathbb N}}

\def\ter{\hfill \vrule width 5 pt height 7 pt depth - 2 pt\smallskip}

\newtheorem{thm}{Theorem}[section]
\newtheorem{corollary}[thm]{Corollary}

\newtheorem{lem}[thm]{Lemma}
\newtheorem{pro}[thm]{Proposition}
\newtheorem{remark}[thm]{Remark}

\begin{document}

\title[Asymptotic convergence to pushed wavefronts ]{Asymptotic convergence to pushed wavefronts in a monostable  equation with delayed reaction}

\author{Abraham Solar}
\address{ Instituto de Matem\'atica y Fisica,
Universidad de Talca, Casilla 747, Talca, Chile}
\ead{asolar.solar@gmail.com}

\author{Sergei Trofimchuk}
\address{ Instituto de Matem\'atica y Fisica,
Universidad de Talca, Casilla 747, Talca, Chile}
\ead{trofimch@inst-mat.utalca.cl}

\begin{abstract} 
We study the asymptotic behavior of solutions to the delayed
monostable  equation $(*)$: $u_{t}(t,x) = u_{xx}(t,x) - u(t,x) + g(u(t-h,x)),$ $x \in \R,\ t >0,$ with monotone reaction term $g: \R_+ \to \R_+$.  Our basic assumption is that this equation possesses  pushed traveling fronts.   First we prove that the pushed wavefronts are 
nonlinearly  stable with asymptotic phase.  Moreover, combinations of these waves attract, uniformly on $\R$, every solution of  equation $(*)$ with the initial datum sufficiently rapidly decaying at one (or at the both)   infinities of 
the real line.  These results  provide a sharp form of the theory of spreading speeds for equation $(*)$. 
  \end{abstract}

\vspace{2pc}

\newpage

\section{Introduction and main results}\vspace{0mm}\noindent
In this work, we study the asymptotic convergence   of  solution $u(t,x)$ of the  initial
value problem for a monostable reaction-diffusion equation with
delayed reaction
\begin{eqnarray} \label{e1}
u_{t}(t,x) &=& u_{xx}(t,x) - u(t,x) + g(u(t-h,x)),
 \\
u(s,x) & =& w_0(s,x), \ s \in [-h,0], \ x \in \R, \label{e2}
\end{eqnarray}
to a combination of traveling waves. 
In the sequel, it is always assumed that the continuous  
function $w_0(s,x)$ is locally H\"older continuous in $x \in \R$,
uniformly with respect to $s$, and  that the function $g:
\mathbb R_+ \to \mathbb R_+$ satisfies the  monostability condition
\vspace{2mm}

\noindent  {\rm \bf(H)}  
the equation $g(x)= x$ has exactly two nonnegative solutions: $0$ and
$\kappa >0$. Moreover, $g$ is $C^1$-smooth in some
$\delta_0$-neighborhood of the equilibria  where $g'(0) >1,$
$g'(\kappa) < 1,$ and also satisfies the Lipshitz condition
$|g(u)-g(v)| \leq L_g|u-v|, \ u,v \in [0,\kappa]$. In addition,
there are $C >0,\ \theta \in (0,1],$ such that   $
\left|g'(u)- g'(0)\right| + |g'(\kappa) - g'(\kappa-u)|\leq Cu^\theta $ for $u\in
(0,\delta_0].$  Without restricting generality, we will also assume that $g$ is linearly and $C^1$-smoothly extended 
on $(-\infty,0]$ and $[\kappa, +\infty)$.  \vspace{3mm}

\noindent   Equation (\ref{e1}) (together with  its non-local versions) is an important model in  the population dynamics \cite{CMYZ,  SEDY, GSW, KGB,  LZh, LLLM, MLLS, MeiI,WLR,YCW,YZ} where it is used to describe  
the spatio-temporal evolution of a single-species population.  In this interpretation of  (\ref{e1}), $g$ is a birthrate function, $u(t,x)$ denotes the population density at  location $x$ and time $t$, and  it is  supposed that the species reaches sexual maturity at age $h > 0$.  Clearly, the Cauchy problem (\ref{e1}), (\ref{e2}) can be solved by  the
method of steps \cite{JYW}, where in the first step we have to look for the
solution of the inhomogeneous linear equation
$$
u_{t}(t,x) = u_{xx}(t,x) - u(t,x) + g(w_0(t-h,x)), \  t \in [0,h],
\ x \in \R,
$$
satisfying  the initial condition $u(0,x) = w_0(0,x)$.  
 Besides the hypothesis {\rm \bf(H)}, from a biological point of view,  it is realistic to assume that the birth function $g$ is either strictly increasing or unimodal (i.e. $g$ has exactly one critical point which is the absolute maximum point \cite{KGB,TT,YCW}) function on $\R_+$.   In the population dynamics,  equation (\ref{e1}) improves 
certain weaknesses (cf. \cite{GBN} or \cite[pp. 56-58]{Tur})  of the logistic growth model given by  the  KPP-Fisher delayed or nonlocal equations  \cite{BNPR, BGHR, DN, FZ, FGT, HT}.  One of the most interesting features of the dynamics in (\ref{e1}) is the existence of smooth positive solutions $u(t,x) = \phi(x+ct)$ satisfying the boundary conditions  $\phi(-\infty)=0$ and $\liminf_{t \to +\infty} \phi(t) >0$ (for $c>0$, cf. \cite{SEDY}). Such solutions are called traveling semi-wavefronts (or wavefronts if additionally $\phi(+\infty) =\kappa$), they describe waves of colonisation propagating with the velocity $c$.  The convergence and stability properties  of  wavefronts to  (\ref{e1}) are quite well understood in the non-delayed case (i.e. when $h=0$).  The studies of the front stability in non-delayed  monostable equation (\ref{e1}) were initiated in 1976 by Sattinger \cite{STG} (see \cite{MeiW} for the state-of-art on this topic), but already the seminal work  of Kolmogorov, Petrovskii, Piskunov (1937)  presented a first deep analysis  of the convergence of the solution $u(t,x)$ of (\ref{e1}), (\ref{e2}) (with $-u+g(u)=u(1-u)$ and  with $w_0(s,x)$ being the Heaviside step function $H(x)$) to a monotone wavefront. 

Now,   the investigation of asymptotic behavior of solution to problem   (\ref{e1}), (\ref{e2})
becomes a much more challenging task when $h >0$. For instance, the recent works \cite{DN,IGT, SEDY,  FGT, HT,LLLM, TT} show that the delay $h$ has a strong influence on the geometry of front's profile $\phi$ and complicates enormously the studies of  the front uniqueness \cite{AGT, CMYZ, TPT, WLR} and stability \cite{CMYZ, KGB,  LZh, LLLM,LvW, MLLS, MeiI,MeiW, WLR}. Moreover, in order to be able to perform  
the local stability analysis of equation (\ref{e1}), it  was always necessary to assume the additional  sub-tangency restriction
\begin{equation} \label{stc}
g(u) \leq g'(0)u, \quad u \geq 0. 
\end{equation}
Under this assumption,  all wavefronts of equation (\ref{e1}) are known as {\it `pulled'} fronts (see \cite{BGHR, ID, RGHK, roth, ES, ST,Xin}  
for  further details),  
model (\ref{e1}) is linearly determined \cite{HR,WLL} and there exists  a positive number  $c_*>0$ (called the minimal speed of propagation) separating the positive axis on the set of admissible  {\it semi-wavefronts}  speeds  $[c_*,+\infty)$  and the set $[0,c_*)$ of velocities $c$ for which does not exist {\it any} non-constant positive bounded wave solution $u(t,x) = \phi(x+ct)$  \cite{SEDY}.   Furthermore,  the minimal speed $c_*$ is determined from the
characteristic equation 
\begin{eqnarray}\label{cra}
\chi(z,c):= z^{2}-cz-1+g'(0)e^{-zch} =0
\end{eqnarray}
as the unique real value $c_\#$ for which  $\chi(z,c)$ has a
positive double zero $\lambda_1(c_\#) = \lambda_2(c_\#)$ (i.e. $c_*$ is equal to
$c_\#$ if (\ref{stc}) holds). Note that for $c > c_\#$ equation $\chi(z,c)=0$
has exactly two positive simple roots, we will denote them as
$\lambda_1(c) < \lambda_2(c)$. 

In this way, as far as we know, all studies  of  wave's stability in
the delayed model (\ref{e1}) have dealt exclusively with the stability of pulled wavefronts. 
Nevertheless, from an ecological point of view, 
models with the birth functions which are not sub-tangential at $u=0$ are also quite interesting  in view  of the
interpretation of non-sub-tangentiality property of $g$ in terms of  a weak Allee effect \cite{BGHR, JP, RGHK}. In the non-delayed case, 
it is well known \cite{ID, RGHK, roth, ES, ST,Xin} that such systems can possess a 
special type of minimal wavefronts  called the {\it `pushed'} fronts.  As the characterising property of a pushed wave for model (\ref{e1}), we can take the following one: the minimal wavefront $u(t,x) = \phi(x+c_*t)$ is pushed if the velocity
$c_*$  is not linearly determined,  i.e. if $c_* > c_\#$. The recent  work \cite{RGHK} explains 
why, contrarily to the pulled waves,  the pushed  colonisation waves  can be considered as waves promoting genetic diversity in the ecological systems.

To the best of our knowledge,  the study of pushed waves in the  monostable 
delayed model (\ref{e1}) was initiated in \cite{LZh, TPT} (curiously, in the first work \cite{sch} dealing with traveling waves in delayed models, all waves were tacitly  presumed to be pulled).  In \cite{TPT},  after assuming  monotonicity of $g$, it was proved that the unique minimal wavefront propagating
with the speed $c_* > c_\#$ must have a strictly increasing profile $\phi$ with  the following asymptotic
representation at $-\infty$:
\begin{equation}\label{af1}
\phi(t+s_0)= e^{\lambda_2 t} + O(e^{(\lambda_2+ \varsigma) t}), \
\lambda_2:= \lambda_2(c_*),\ \varsigma >0, \quad t \to -\infty.
\end{equation}
It should be noted that the situation 
when non-monotone (for example, unimodal)  birth function $g: [0,\kappa] \to \R_+ $ does not satisfy (\ref{stc}) is not  completely understood till now. In fact, even the existence of the minimal speed of propagation $c_*$, as the lowest value from  a closed connected unbounded set of all admissible wavefront (or semi-wavefront \cite{BNPR, SEDY}) velocities,
is not  yet proved for the case of  non-monotone and not sub-tangential $g$.  From the formal point of view, the existence of the pushed fronts to the delayed model (\ref{e1}) neither was established  in \cite{TPT}. In any case,  this point can be easily completed:  
\begin{pro}\label{mainP} Assume that $u = \phi(x +c_*t), \ c_* > c_\#(h_0)$,  is a pushed traveling front  to the monotone model (\ref{e1}) considered with some fixed $h_0\geq 0$.  Then there exists a positive $\delta$ such that  equation (\ref{e1}) possesses a pushed traveling front for each non-negative $h \in (h_0-\delta, h_0+\delta)$.  In particular,  there exists a delayed equation (\ref{e1}) with $h >0$ possessing the minimal monotone wavefront
$u = \phi(x +c_*t)$ with the profile $\phi$ satisfying the asymptotic
formula (\ref{af1}).
 \end{pro}
 {\it Proof. }
Since $c_\#(h)$ depends continuously on $h \geq 0$, the first part
of Proposition \ref{mainP} will be proved if we establish the
lower semicontinuity of $c_*(h)$ at $h_0$.  Then the existence of
pushed wavefronts to the  equation (\ref{e1}) considered with
small positive delays follows from the existence of the pushed
wavefronts to  the Fisher type population genetic model
\cite[Theorem 11]{HR}
$
u_{t}(t,x) = u_{xx}(t,x) - u(t,x) +
(10u(t,x)+3u^2(t,x)-5u^3(t,x))/{8}.
$
Hence, it suffices to prove the following

\vspace{2mm}

\noindent {\bf Claim.} {\it Suppose that $h_j \to h_0, \ c_*(h_j) \to c_0$ as $j \to +\infty$. Then $c_0 \geq c_*(h_0)$.}

\vspace{2mm}

\noindent Indeed, take some $c > c_0$. Then, for all sufficiently large $j$, the equation
$$
u_{t}(t,x) = u_{xx}(t,x) - u(t,x) + g(u(t-h_j,x))
$$
has a unique (up to translation) positive strictly monotone
wavefront $u(t,x) = \phi_j(x+ct)$. Without the loss of the
generality, we can assume that $\phi_j(0)=\kappa/2$. It is easy to
see (cf. \cite{TPT}) that each profile $\phi_j$ satisfies the
integral equation
\begin{equation}\label{ie1} \hspace{-2cm}
\phi(t) = \frac{1}{\xi_2 - \xi_1} \left(\int^t_{-\infty} e^{\xi_1
(t-s)}g(\phi(s-ch_j))ds + \int_t^{+\infty}e^{\xi_2
(t-s)}g(\phi(s-ch_j))ds\right),
\end{equation}
where $\xi_1< 0<\xi_2$ are roots  of the  equation $z^2 -c z -1
=0$.  Since $|\phi'_j(t)| \leq \kappa/\sqrt{c^2+4}$, $|\phi_j(t)|
\leq \kappa$, the sequence $\phi_j$ has a subsequence $\phi_{j_k}$
wich converges, uniformly on compact subsets of $\R$, to the
monotone continuous bounded function $\phi_0(t), \ \phi_0(0) =
\kappa/2$.  By the Lebesgue«s dominated convergence theorem,
$\phi_0$ satisfies the equation (\ref{ie1}) with $h_0$ and
therefore  $\phi_0$ is a positive profile of strictly 
monotone wavefront propagating with the velocity $c$ \cite{SEDY, TPT}. In this way,
$c \geq c_*(h_0)$ for every $c> c_0$ that yields $c_0 \geq
c_*(h_0)$.
\hfill \ter 

\vspace{2mm}

Formula (\ref{af1}) implies that pushed profiles $\phi(s)$
converges to $0$ at $ -\infty$ more rapidly than the
profiles of other (i.e. non-minimal or pulled) waves behaving as 
$$
\phi(t+s_0)= (-t)^me^{\lambda_1 t} + O(e^{(\lambda_1+ \varsigma) t}), \
\lambda_1:= \lambda_1(c), \ \varsigma >0, \ m \in \{0,1\}, \quad t \to -\infty.
$$
The fast asymptotic decay of pushed fronts at $-\infty$ makes them similar to the so-called {\it bistable} fronts \cite{FML, SZ, Xin}. Actually, by analysing the inside dynamics of wavefronts, Garnier {\it et al} \cite{ID} (in the non-delayed case) and  Bonnefon {\it et al} \cite{BGHR} (in the delayed case) have recently proposed  a general definition of pushed waves which allows to consider the monostable pushed fronts and the bistable fronts within a unified framework.   An additional argument if favor of this  insight is provided by the theory of nonlinear stability of waves. Indeed, both monostable pushed fronts and bistable fronts are proved to have rather good stability properties \cite{FML,OM, roth, SZ, ST, ST2}. Furthermore,  the most complete and comprehensible proof 
of the asymptotic stability of monostable pushed front given in \cite{roth} uses constructions and results obtained  for a bistable model in \cite{FML}.  

Hence, the main aim of the present paper is to study the stability properties of monostable pushed fronts 
to the monotone delayed model (\ref{e1}).  We are going to achieve this goal 
by developing several ideas and methods from \cite{FML, OM, roth, TPT}. We also will establish the asymptotic 
convergence of solutions for the initial value problem   (\ref{e1}), (\ref{e2}) to an appropriate pushed wavefront 
when, in addition to {\rm \bf(H)},   $g$ is monotone   and  when $w_0$ satisfies,  for some
$A, B >0, \ \sigma \in (0,\kappa)$  and $\mu > \lambda_1(c_*)$  the following  conditions {\rm \bf(IC)}:
\begin{description}
 \item [$(IC1)$] $\quad 0\leq w_0(s,x)\leq\kappa, \quad  x \in
\mathbb{R}, \quad s\in[-h,0];$ \item[$(IC2)$]  $\quad w_0(s,x)\leq
Ae^{\mu x}, \quad x\in\mathbb{R}, \quad s\in[-h,0];$ \item[$(IC3)$]
$\quad w_0(s,x)>\kappa - \sigma, \quad s\in[-h,0], \ \ x \geq B$.\end{description}
From the monotonicity of $g$  and the hypotheses  {\rm \bf(H)}, {\rm \bf(IC)}, by  invoking the well-known existence and uniqueness results
and the comparison principle \cite[Chapter 1, Theorems 12,
16]{AF},  we can deduce  the existence of a unique classical solution
$u= u(t,x): [-h,+\infty) \times \R \to [0,\kappa]$ to (\ref{e1}),
(\ref{e2}) (i.e. of a continuous bounded function $u$ having continuous
derivatives $u_t, u_x, u_{xx}$ in $\Omega = (0,+\infty)\times \R$
and satisfying (\ref{e1}) in $\Omega$ as well as (\ref{e2}) in
$[-h,0]\times \R$). As the following proposition shows, the asymptotic behavior of this solution $u(t,x)$ on bounded subsets of $\R$ 
is quite simple: 
\begin{pro}\label{mainP2}  Suppose that  the initial datum $w_{0}\not\equiv 0$  satisfies $(IC1)$ and that the Lipshitz continuous map $g:[0,\kappa]\to [0,\kappa]$ has exactly two fixed points: $0$ and
$\kappa >0$. Then
$
\lim_{t\rightarrow\infty}u(t,x)=\kappa
$
uniformly on compact subsets of $\R$.
\end{pro}

At first glance, if additionally we assume the monotonicity of $g$, Proposition \ref{mainP2} seems to follow from quite general results on spreading speeds to 
continuous-time semiflows established in  \cite{LZhA, LZh}. Indeed, \cite[Theorem 34]{LZh} shows that  even 
rather weak positivity condition  assumed 
in Proposition \ref{mainP2} is enough to assure stronger convergence 
\begin{equation}\label{sps}
\lim_{t\to\infty}\sup_{x \in [-c't,c't]}|u(t,x)-\kappa|=0, \ c' \in (0,c_*), 
\end{equation}
once $g$ is a subhomogeneous function: $\rho g(x) \leq  g(\rho x) $ for all $\rho \in [0,1]$  and $x \geq 0$.  It is easy to see, however, that the latter condition implies  the sub-tangency inequality (\ref{stc}). 

Our proof of Proposition \ref{mainP2} follows closely the main lines of  \cite{AW}, where 
Aronson and Weinberger  established a similar result for non-delayed equations.  
See also \cite[Theorem 3.2]{YZ} for an analogous assertion proved  
 for a non-diffusive delay differential equation with spatial non-locality in 
an unbounded domain.  In general (e.g. under condition $(IC2)$) the convergence of $u(t,\cdot)\to \kappa$, $t\to +\infty$,  is not uniform on $\R$: this  is an immediate outcome of  our subsequent investigation of  the 
asymptotic behavior of   the {\it entire solution} $u(t,x)$ as $t\to +\infty$ on {\it  the whole real $x$-line $\R$. }

In order to state the main  results of this paper, we take a  pushed front $\phi(x+c_*t)$ for equation  (\ref{e1})  and fix a positive number $\lambda < \mu$ such
that $\lambda \in (\lambda_1(c_*), \lambda_2(c_*))$. We will  also consider the  Banach space
\[
C_{\lambda}(\R) = \left\{y \in C(\R, \R): |y|_{\lambda}: =
\max\{\sup_{x \leq 0}e^{-\lambda x} |y(x)|, \sup_{x \geq 0}
|y(x)|\} < \infty  \right\}.
\]
Observe that $|y|_\lambda = \sup_{x \in \R} |y(x)|/\eta(x)$, where
$\eta(x):=\min\{  e^{\lambda x}, 1\}$. 
Our first theorem shows that the pushed front $\phi(x +c_*t), \ c_* > c_\#,$ is nonlinearly  stable with asymptotic phase \cite{BS}:
\begin{thm} \label{cocoT} Let $g$ be monotone and  conditions \mbox{\rm ({\bf IC})}, \mbox{\rm ({\bf H})}  be satisfied. Then for every $\epsilon >0$ there exists $\delta >0$ such that $|\phi(\cdot+c_*s) - w_0(s,\cdot)|_\lambda < \delta,$ $s \in [-h,0],$ implies that $|\phi(\cdot +c_*t) - u(t,\cdot)|_\lambda < \epsilon$ for all $t \geq 0$. Here $u(t,x)$ is  solution  of the initial value problem (\ref{e1}), (\ref{e2}). Furthermore, there exists $s_0$ such that  $|\phi(\cdot +c_*t+s_0) - u(t,\cdot)|_\lambda \to 0$ as $t \to +\infty$. 
\end{thm}
The stability result of Theorem \ref{cocoT} follows from Corollary \ref{coco} proved in Section 2 while the  asymptotic convergence $u(t,x) \to  \phi(x +c_*t+s_0), \ t \to +\infty$, follows from  the next  theorem. It describes the global stability properties of the pushed fronts with respect to initial data satisfying the hypothesis   ({\bf IC}): 
\begin{thm}\label{main1} Let $g$ be monotone and conditions \mbox{\rm ({\bf IC})}, \mbox{\rm ({\bf H})}  be satisfied. Then  the solution $u(t,x)$ of the initial value problem (\ref{e1}), (\ref{e2}) asymptotically converges to  a shifted front. In fact,  for some $s_0\in \R$,
\begin{equation}\label{mcr}
\lim_{t \to \infty }\sup_{x\in \R}|u(t,x) -   \phi(x +c_*t+ s_0)|/\eta(x+c_*t) =0.
\end{equation}
 \end{thm}
It is instructive to compare Theorems \ref{cocoT},  \ref{main1} with stability results obtained for non-critical 
pulled fronts in the delayed model (\ref{e1}) with monotone reaction $g$ satisfying (\ref{stc}) and {\rm \bf(H)}. For example, taking initial functions $w_0$ satisfying  $(IC1)$ and assuming that the initial disturbance
$\phi(\cdot+cs) - w_0(s,\cdot)$ belongs to the weighted Sobolev space  $H^1_{\eta^2}(\R)$ and depends continuously on $s\in [-h,0]$, 
Mei et al \cite[Theorem 2.2]{MLLS} proved  that $|\phi(\cdot +ct) - u(t,\cdot)|_0 \to 0$ exponentially when $t\to +\infty$.  Hence, in view of the continuous imbedding $H^1_{\eta^2}(\R) \subset C_{\lambda}(\R)\cap C^{0,1/2}(\R_+)$,  initial functions $w_0(s,x)$ in \cite{MLLS} are uniformly H\"older continuous in $x$ and converge  at $+\infty$, $w_0(s,+\infty) = \kappa$ (in fact, this convergence is uniform in $s \in [-h,0]$,  so that each $w_0$ meets trivially the restriction (IC3)). 
They should also satisfy the inequality

\begin{equation}\label{icc}\hspace{-2.5cm}
|\phi(x+cs) - w_0(s,x)| \leq Ce^{\lambda x}, \ x \in \R,\ s \in [-h,0], \ \mbox{for some} \ C>0, \ \lambda \in (\lambda_1, \lambda_2).
\end{equation}
Due to the asymptotic representation (\ref{af1}) and to certain freedom in the choice of $\lambda, \mu$, in the case of pushed fronts, the latter condition amounts precisely to the hypothesis (IC2).  Nevertheless, in  contrast to inequality  (\ref{icc}) considered with a pushed front $u= \phi(x+c_*s)$, the same inequality  considered with a  pulled front $u= \phi(x+cs)$
 is not satisfied if we take the Heaviside step function $H(x)$ as the initial function $w_0(s,x)= H(x)$. Thus the question about the asymptotic form
 of solution $u(t,x)$ to the Cauchy problem  (\ref{e1}), (\ref{e2}) with $w_0(s,x)= H(x)$ and with the sub-tangential $g$ still remains unanswered in the delayed case.  It is worth to recall that precisely this question formulated for a non-delayed monostable equation (\ref{e1}) was the main object of studies in the seminal work by Kolmogorov, Petrovskii, Piskunov in 1937.  
 
Now, it is worth noticing that equation (\ref{e1}) is invariant wit respect to the transformation $x \to -x$  so that 
the statements of Theorems \ref{cocoT} and \ref{main1} can be easily adapted to the case when the initial function $w_0(s,-x)$  meets the hypothesis ({\bf IC}). Evidently, in such a case,  we should use {\it pushed backs} of the form $u =  \phi(-x +c_*t)$ instead of the pushed wavefronts.  Then the natural question is whether solution $u(t,x)$  converges to a combination of a pushed  front and a pushed back  when  the both non-zero functions $w_0(s,x), w_0(s, -x)$ satisfy  conditions (IC1), (IC2). In particular, this happens when  $w_0$ has  compact support. 
To the best of our
knowledge, the studies of the asymptotic form of solutions to the  monostable reaction-diffusion  equations  having compactly supported initial data were initiated  in \cite{AW,roth, ST2,UC}. Here, we analyse a similar problem  in the presence of delay; hence, our third theorem considers the initial data for  (\ref{e1}), (\ref{e2})  exponentially vanishing at both  infinities.
\begin{thm}\label{main2} Assume that $u = \phi(x +c_*t), \ c_* > c_\#$,  is a pushed traveling front  to equation (\ref{e1}). If non-zero functions $w_0(s,x), w_0(s, -x)$ satisfy  conditions (IC1), (IC2) then  the solution $u= u(t,x)$ of the initial value problem (\ref{e1}), (\ref{e2}) asymptotically converges to  a  combination of two shifted fronts, i.e. for some $s_1, s_2 \in \R$,
\begin{eqnarray*} 
\lim_{t \to \infty }\sup_{x\leq 0} |u(t,x) -
\phi(x +c_*t+ s_1)|/\eta(x+c_*t)=0,\\
\lim_{t \to \infty }\sup_{x\geq
0}|u(t,x) -   \phi(-x +c_*t+ s_2)|/\eta(-x+c_*t) =0.
\end{eqnarray*}
\end{thm}

Clearly, Theorem \ref{main2} combined with the comparison principle shows that relation (\ref{sps}) holds for each solution $u= u(t,x)$ to (\ref{e1}) once associated initial datum $w_0(s,x)\not\equiv 0$ satisfies  (IC1). Moreover,  since Theorem \ref{main2}  implies that 
$$
\lim_{t\to\infty}\sup_{x \not\in (-c't,c't)}u(t,x)=0, \ c' >c_*, 
$$
we can conclude that  the speed $c_*$ of pushed waves coincides with the spreading speed for model (\ref{e1}). Without restriction (\ref{stc}), this  important result  was for the first time  established in \cite{LZhA, LZh}  (in a much more general setting).  Therefore Theorem \ref{main2}  can be also viewed as an essential  improvement of the mentioned Liang and Zhao result for the particular case of Eq. (\ref{e1}). 

As in \cite{FML, roth}, the method of sub- and super-solutions is 
a key tool for proving our main results. The sub- and super-solutions will be obtained 
as suitable deformations (invented by Fife and McLeod in \cite{FML} for the bistable systems and  adapted by Rothe in \cite{roth} for the monostable  equations) of the pushed wavefront.  The other important idea exploited in    
\cite{FML, roth} is  the use of an appropriate Lyapunov functional for  proving the wave stability.  However, the construction of such a functional seems to be a rather difficult task in the case of the functional differential equation (\ref{e1}).  Thus, instead of this, we decided to use the Berestycki and Nirenberg method of the
sliding solutions \cite{BN,TPT} as well as some ideas of the  approach
developed by Ogiwara and Matano in \cite{OM}. It is natural to expect that the  rate of convergence 
in (\ref{mcr}) is exponential, see  e.g.  \cite{FML, MLLS, MeiI, roth, BS}.  The demonstration of this fact, however, 
is based on a different approach and will be considered  in a separate work.

Finally, we say a few words about the organization of the paper. The results of 
Theorems \ref{cocoT} and   \ref{main1} follow from Corollary \ref{coco} and Theorem \ref{cocos} which are proved in Section \ref{sub1}. Then various auxiliary 
results are proved in Section  \ref{sub2}  (Proposition \ref{mainP2}) and Section  \ref{sub3} (an important stability 
Lemma \ref{lem4} among others).  In the last section of the paper, we completes the proof of Theorem \ref{main2}.

\section{Proof of Theorems \ref{cocoT}  and \ref{main1}} \label{sub1}
Let  $u = \phi(x+c_*t), \ c_* > c_\#,$ be a pushed traveling
front  to equation (\ref{e1}). In the sequel, to simplify  the
notation, we will avoid  the subscript  $*$ in $c_*$ so that
$u= \phi(x+c_*t) = \phi(x+ct)$.  As it is usual, we consider
the moving coordinate frame $(t,z)$ where  $z= x+ct$. Set $w(t,z)
= u(t,z - ct)$,  then equation (\ref{e1})  takes the form
\begin{eqnarray}\label{wz}
& & w_{t}(t,z)= w_{zz}(t,z)-cw_{z}(t,z)-w(t,z)+g(w(t-h, z-ch)), \\
& & w(s,z)  = \tilde w_0(s,z):= w_0(s,z-cs), \ s \in [-h,0], \ z \in \R.
\label{wzic}
\end{eqnarray}
First, following Fife and McLeod \cite[Lemma 4.1]{FML} and
Rothe \cite[Lemma 1]{roth},  we  prove the next assertion.
\begin{lem} \label{lem1} Assume that the hypothesis {\rm ({\bf H})} is satisfied.
Then there exist positive constants $\gamma, C, q_0^+$ (depending
only on $g, \phi, c, h, \lambda$) and $q_0^- = \sigma$ such that
the inequality
\begin{equation}\label{mlem}0 \leq w(s,z)\leq \phi(z)+q \eta(z), \quad z\in\mathbb{R},\quad s\in [-h,0],
\end{equation}
with $q \in (0,q_0^+]$ implies
 \begin{equation}\label{mlemm}0 \leq w(t,z)\leq \phi(z+Cq)+qe^{-\gamma t}
\eta(z), \quad z\in\mathbb{R},\quad t\geq -h.
\end{equation}
Similarly, the inequality
\begin{equation}\label{mlemN}\phi(z+Cq)-q \eta(z) \leq w(s,z)\leq \kappa,  \quad z\in\mathbb{R},\quad s\in [-h,0],
\end{equation}
with $q \in (0,q_0^-]$ implies
 \begin{equation}\label{mlemmN}
 \phi(z)-qe^{-\gamma t} \eta(z) \leq w(t,z)\leq \kappa,
\quad z\in\mathbb{R},\quad t\geq -h.
\end{equation}
\end{lem}
{\it Proof. }
For the
convenience of the reader, the proof is divided into five steps. Recall that the positive numbers $\delta_0, \sigma$ are defined in 
({\bf H}) and ({\bf IC}), respectively.   
\vspace{2mm}

\noindent \underline{Step I}. We claim that given $\sigma \in
(0,\kappa)$, there are positive $\delta_1^*< \delta_0,\
\gamma_1^*< \lambda c$ such that
$$
g(u)- g(u- qe^{\gamma h}) \leq q(1-2\gamma), \ \mbox{for all} \
(u,q)\in[\kappa-\delta_1^*,\kappa + \delta_1^*]\times[0,\sigma], \
\gamma \in [0,\gamma_1^*].
$$
Indeed, it suffices to note that, given $\sigma \in (0,\kappa)$,
the continuous function
\begin{eqnarray*} \hspace{-2.5cm}&&  G(u,q, \gamma):=\left\{
\begin{array}{ll} 1+(
g(u- e^{\gamma h}q)-g(u))/q, & {(u,q)\in[\kappa-\delta_1^*,\kappa+ \delta_1^*]\times (0, \sigma]}, \ \gamma \in [0,\gamma_1^*]; \\
    1- e^{\gamma h }g'(u), & {u\in[\kappa-\delta_1^*,\kappa+ \delta_1^*], \ q =0},\ \gamma \in [0,\gamma_1^*],
\end{array}%
\right.
\end{eqnarray*}
satisfies  $G(\kappa, q, 0) > 2 \gamma_1^*, \ q \in [0,
\sigma],$  for sufficiently small $\gamma_1^*, \delta_1^*$ (recall
that $g'(\kappa)<1$). Thus  $G(u,q, \gamma) > 2\gamma$ for all 
$(u,q)\in[\kappa-\delta_1^*,\kappa + \delta_1^*]\times[0,\sigma], \
\gamma \in [0,\gamma_1^*]$ if $\gamma_1^*, \delta_1^*$  are sufficiently small.

\vspace{2mm}

\noindent \underline{Step II}.  As in \cite{FML,roth}, we have to construct appropriate super- and sub-solutions.  Consider the
nonlinear operator ${\mathcal N}$ defined as
\begin{equation*} \hspace{-1.5cm} 
{\mathcal N}w(t,z):= w_{t}(t,z)
-w_{zz}(t,z)+cw_{z}(t,z)+w(t,z)-g(w(t-h, z-ch)). 
\end{equation*}
By definition,  continuous function $w_+:  \R_+\times \R \to \R$ is
called a super-solution for (\ref{wz}), if, for some $z_* \in \R$,
this function  is $C^{1,2}$-smooth in the domains $\R_+\times
(-\infty, z_*]$ and $\R_+\times [z_*, +\infty)$ and
\begin{equation}  \label{sso} \hspace{-2.5cm} {\mathcal N}w_+(t,z) \geq 0 \ 
\mbox{for }  t>0,  \ z \not= z_*,  \ \mbox{while} \
(w_+)_z(t,z_*-)\geq (w_+)_z(t,z_*+)\ \mbox{for }\ t >0.
\end{equation}
Sub-solutions $w_-$ are defined analogously, with the inequalities
"$\geq$" reversed in (\ref{sso}). We will look for super- and sub-solutions of the form
\begin{equation*} \hspace{-1.5cm}
w_{+}(t,z):= \phi(z+\epsilon(t)) + qe^{-\gamma t}\eta(z), \quad
w_{-}(t,z):= \phi(z-\epsilon_1(t)) - qe^{-\gamma t}\eta(z),
\end{equation*}
where, for appropriate positive parameters $\alpha,\gamma$ (to be fixed
later and depending only on $g, \phi, c, h, \lambda$), increasing
$\epsilon(t),\ \epsilon_1(t)$ are defined by
\begin{eqnarray}\hspace{-2.5cm}
\epsilon(t):= \frac{\alpha q}{\gamma}(e^{\gamma h}- e^{-\gamma
t})>0,  \quad \epsilon_1(t): = \epsilon(t)-\epsilon(+\infty) =
-\frac{\alpha q}{\gamma}e^{-\gamma t}<0, \quad  t > -h. \nonumber
\end{eqnarray}
Note that the smoothness conditions and  the second inequality in
(\ref{sso}) with $z_* =0$ are obviously fulfilled because of
$$\frac{\partial w_+(t,0+)}{\partial z}- \frac{\partial
w_+(t,0-)}{\partial z}= -  q\lambda e^{-\gamma t} <0, \quad
\frac{\partial w_-(t,0+)}{\partial z}- \frac{\partial
w_-(t,0-)}{\partial z}=  q\lambda e^{-\gamma t} >0,$$ so that we
have to check the first inequality of (\ref{sso}) only.
Since $g, \phi, \epsilon$ are strictly increasing, we have, for $z
\not=0$, that
\begin{eqnarray} \hspace{-2.5cm}
{\mathcal N}w_+(t,z):= \epsilon'(t)\phi'(z+\epsilon(t))- \gamma q
e^{-\gamma t}\eta(z)-\phi''(z+\epsilon(t))- q e^{-\gamma
t}\eta''(z)+c\phi'(z+\epsilon(t)) \nonumber \\ \hspace{-2.5cm}
+cq e^{-\gamma
t}\eta'(z)+\phi(z+\epsilon(t))+ q e^{-\gamma t}\eta(z)-g(w_{+}(t-h,z-ch))
\geq \alpha  q e^{-\gamma t}\phi'(z+\epsilon(t))\nonumber\\ \hspace{-2.5cm}
 - \gamma q e^{-\gamma t}\eta(z)+ cq e^{-\gamma t}\eta'(z) +q e^{-\gamma t}\eta(z)-
q e^{-\gamma
t}\eta''(z)\nonumber\\ \hspace{-2.5cm}
+g(\phi(z-ch+\epsilon(t)))-g(\phi(z-ch+\epsilon(t)) +
qe^{-\gamma (t-h)}\eta(z-ch)) \nonumber;
\end{eqnarray}
\begin{eqnarray} \hspace{-2.5cm}
{\mathcal N}w_-(t,z):= -\epsilon_1'(t)\phi'(z-\epsilon_1(t))+
\gamma q e^{-\gamma t}\eta(z)-\phi''(z-\epsilon_1(t)) +q
e^{-\gamma t}\eta''(z)+c\phi'(z-\epsilon_1(t)) \nonumber \\ \hspace{-2.5cm}
-cq e^{-\gamma
t}\eta'(z)+\phi(z-\epsilon_1(t))- q e^{-\gamma t}\eta(z)-g(w_{-}(t-h,z-ch))
\leq -\alpha  q e^{-\gamma t}\phi'(z-\epsilon_1(t))\nonumber\\ \hspace{-2.5cm}
+\gamma q e^{-\gamma t}\eta(z)- cq e^{-\gamma t}\eta'(z) -  q e^{-\gamma t}\eta(z)+
q e^{-\gamma
t}\eta''(z)+g(\phi(z-ch-\epsilon_1(t)))\nonumber
 \\ \hspace{-2.5cm} -g(\phi(z-ch-\epsilon_1(t))
- qe^{-\gamma (t-h)}\eta(z-ch)).\nonumber
\end{eqnarray}
%
Since $\lambda \in (\lambda_1(c), \lambda_2(c))$ and $g'(0)
>1$, we can choose  sufficiently small $\gamma \in  (0,\gamma_1^*)$ and $\delta \in (0, \kappa/2)\cap
(0,\delta_1^*)\cap (0,\sigma)$,  such that, for all $\bar s <\delta$ it holds
\begin{eqnarray}\label{C1a}
-\lambda^2+c\lambda +1-\gamma-g'(\bar s)e^{-\lambda ch+\gamma
h}>0.
\end{eqnarray}
In addition, we can take $\delta$ such that the unique real roots
$z_0< z_1  < z_2$ of the equations
$$
\phi(z_0) = \delta/4, \quad \phi(z_1) + 0.25 \delta
\eta(z_1)=\delta/2; \quad \phi(z_2) = \kappa -\delta/2,
$$
are such that  $z_1 <-ch < 0< z_2$.  
From now on, we will fix $\alpha, q_0^\pm$ defined by
$$q_0^+ = \delta e^{-\gamma h}/2, \quad q_0^- = \sigma,  \quad
\alpha= (\gamma+ e^{\gamma h}L_g)/\beta, \quad \mbox{with} \
\beta:= \min_{z\in[z_{0},z_{2}+ch]}  \phi'(z).
$$
We observe that $\alpha, q_0^\pm$ and $\gamma$ depends only on $g,
\phi, c, h, \lambda, \sigma$.

\vspace{2mm}

\noindent \underline{Step III}. We claim that  ${\mathcal
N}w_+(t,z) \geq 0$ for all $z\not=0, \ t \geq 0$ and $q \leq q_0^+=
\delta e^{-\gamma h}/2$.

Indeed, suppose first that  $z-ch+\epsilon(t) \leq z_1$, then $z
\leq z_1+ch- \epsilon(t) < -\epsilon(t) <0$ and
$$
\phi(z-ch+\epsilon(t)) < \delta/2,  \ \phi(z-ch+\epsilon(t)) +
qe^{-\gamma (t-h)}\eta(z-ch) < \delta.
$$
As a consequence, we can invoke the mean value theorem and and (\ref{C1a}) to conclude
that, for some $\bar s \in (0,\delta)$,
\begin{eqnarray*}
&& \hspace{-2.5cm}g(\phi(z-ch+\epsilon(t)))-g(\phi(z-ch+\epsilon(t)) +
qe^{-\gamma (t-h)}\eta(z-ch)) = - qe^{-\gamma
(t-h)}\eta(z-ch)g'(\bar s),\\ && \hspace{-2.5cm} {\mathcal N}w_+(t,z) \geq -\gamma
q e^{-\gamma t}\eta(z)+ cq e^{-\gamma t}\eta'(z) +
 q e^{-\gamma t}\eta(z)-
q e^{-\gamma
t}\eta''(z)- qe^{-\gamma (t-h)}\eta(z-ch)g'(\bar s) \\
&& \hspace{-2.5cm} =q e^{-\gamma t}\left( [1- \gamma] \eta(z)+ c\eta'(z)
-\eta''(z)- e^{\gamma h}\eta(z-ch)g'(\bar s)\right) \\
&& \hspace{-2.5cm} = q e^{-\gamma t + \lambda z}\left( 1- \gamma + c\lambda
-\lambda^2- e^{\gamma h-\lambda ch}g'(\bar s)\right) >0.
\end{eqnarray*}
Similarly, if $z-ch+\epsilon(t) \geq z_2,$ then we have that
$$
 \kappa+\delta/2 \geq \phi(z-ch+\epsilon(t)) + qe^{-\gamma (t-h)}\eta(z-ch) \geq \phi(z-ch+\epsilon(t))\geq  \kappa- \delta/2.
$$ Therefore, due to  Step I and (\ref{C1a}),  for all $t \geq 0$,
$$  g(\phi(z-ch+\epsilon(t)))-g(\phi(z-ch+\epsilon(t)) + qe^{-\gamma (t-h)}\eta(z-ch)) \geq -qe^{-\gamma t}\eta(z-ch)(1-2\gamma),$$
\begin{eqnarray*}  &&\hspace{-2.5cm}
{\mathcal N}w_+(t,z) \geq q e^{-\gamma t}\left( [1- \gamma]
\eta(z)+ c\eta'(z)
-\eta''(z)- (1-2\gamma)\eta(z-ch)\right) \geq \\
&&\hspace{-2.5cm}  q e^{-\gamma t} \left\{
\begin{array}{ll} e^{\lambda z}
     [1- \gamma+ c\lambda
-\lambda^2- e^{-\lambda ch}(1-2\gamma)], & {z< 0} \\
     \gamma, & {z > 0}
\end{array}%
\right\} >0.
\end{eqnarray*}
Finally, if $z_{1}< z-ch+\epsilon(t) <  z_{2}$, we find that
$$
\phi(z-ch+\epsilon(t)) < \kappa- \delta/2,  \quad
\phi(z-ch+\epsilon(t)) + 0.25 \delta \eta(z-ch+\epsilon(t))
>\delta/2$$ so that $\phi(z-ch+\epsilon(t))  >\delta/2 - 0.25
\delta \eta(z-ch+\epsilon(t)) > \delta/4$, and $ z+\epsilon(t) \in
[z_0,z_2+ch] $. Obviously,
$$
|g(\phi(z-ch+\epsilon(t)))-g(\phi(z-ch+\epsilon(t)) + qe^{-\gamma
(t-h)}\eta(z-ch))| \leq L_gqe^{-\gamma (t-h)}\eta(z-ch).
$$
Therefore, since $\eta(z) + c\eta'(z) -\eta''(z) >0 $ for $z \not=0$ and $\eta(z)
\in (0,1],$ we get
\begin{eqnarray} \hspace{-.5cm}
{\mathcal N}w_+(t,z) \geq q e^{-\gamma t}\left\{\alpha\beta +
\eta(z) + c\eta'(z) -\eta''(z) - \gamma -e^{\gamma h}L_g\right\}> 0
\nonumber
\end{eqnarray}

\vspace{1mm}

\noindent for all $t\geq 0$.
Hence, there exist some constants  $\alpha, \gamma, q_0^+ >0,$
depending only on the wavefront profile $\phi$, the nonlinearity
$g$ and $c,h, \lambda$ such that,  for any choice of  $q \in
(0,q_0^+)$ it holds  ${\mathcal N}w_+(t,z)>0$  for all $t\geq 0$ and
$z\not=0$. This proves the first inequality in (\ref{sso}).

\vspace{2mm}

\noindent \underline{Step IV}. We claim that  ${\mathcal
N}w_-(t,z) \leq 0$ for all $z\not=0, \ t \geq 0$ and $q \leq
q_0^-= \sigma$.

\noindent 
Indeed, suppose first that  $z-ch-\epsilon_1(t) \leq z_1$, then
$z \leq z_1 +\epsilon_1(t) +ch <0 $ and 
$$\phi(z-ch-\epsilon_1(t)) - qe^{-\gamma (t-h)}\eta(z-ch)
<\phi(z-ch-\epsilon_1(t)) < \delta/2.
$$
As a consequence,  the mean value theorem yields 
that for some $\bar s  <\delta$,
\begin{eqnarray*}
&&  \hspace{-2.5cm} g(\phi(z-ch-\epsilon_1(t)))-g(\phi(z-ch-\epsilon_1(t)) -
qe^{-\gamma (t-h)}\eta(z-ch)) =  qe^{-\gamma
(t-h)}\eta(z-ch)g'(\bar s),\\ &&  \hspace{-2.5cm} {\mathcal N}w_-(t,z) \leq \gamma
q e^{-\gamma t}\eta(z)-cq e^{-\gamma t}\eta'(z) -
 q e^{-\gamma t}\eta(z)+
q e^{-\gamma
t}\eta''(z)+ qe^{-\gamma (t-h)}\eta(z-ch)g'(\bar s) \\
&&  \hspace{-2.5cm} =-q e^{-\gamma t}\left( [1- \gamma] \eta(z)+ c\eta'(z)
-\eta''(z)- e^{\gamma h}\eta(z-ch)g'(\bar s)\right)  \\
&&  \hspace{-2.5cm} = -q e^{-\gamma t} e^{\lambda z}\left( 1- \gamma + c\lambda
-\lambda^2- e^{\gamma h-\lambda ch}g'(\bar s)\right) <0.
\end{eqnarray*}
Similarly, if $z-ch-\epsilon_1(t) \geq z_2,$ then  we have that
$\phi(z-ch-\epsilon_1(t))\geq  \kappa- \delta/2$ and therefore 
\begin{equation*} \label{C1b} \hspace{-2.5cm}
g(\phi(z-ch-\epsilon_1(t)))-g(\phi(z-ch-\epsilon_1(t)) -
qe^{-\gamma (t-h)}\eta(z-ch)) \leq (1-2 \gamma) qe^{-\gamma
t}\eta(z-ch)
\end{equation*}
for all $t \geq 0, \ q \in [0,\sigma]$.
In consequence, 
\begin{eqnarray*} &&
{\mathcal N}w_-(t,z) \leq -q e^{-\gamma t}\left( [1- \gamma]
\eta(z)+ c\eta'(z)
-\eta''(z)- (1- 2\gamma)\eta(z-ch)\right) \leq \\
&& -q e^{-\gamma t} \left\{
\begin{array}{ll} e^{\lambda z}
     [1- \gamma+ c\lambda
-\lambda^2- e^{-\lambda ch}(1-2\gamma)], & {z< 0} \\
     \gamma, & {z >0 }
\end{array}%
\right\} <0.
\end{eqnarray*}
Finally, if $z_{1}< z-ch-\epsilon_1(t) <  z_{2}$, we find that
$$
\phi(z-ch-\epsilon_1(t)) < \kappa- \delta/2,  \quad
\phi(z-ch-\epsilon_1(t)) + 0.25 \delta \eta(z-ch-\epsilon_1(t))
>\delta/2$$ so that $\phi(z-ch-\epsilon_1(t))  >\delta/2 - 0.25
\delta \eta(z-ch-\epsilon_1(t)) > \delta/4$, and $ z-\epsilon_1(t)
\in [z_0,z_2+ch] $. Obviously,
$$
|g(\phi(z-ch-\epsilon_1(t)))-g(\phi(z-ch-\epsilon_1(t)) -
qe^{-\gamma (t-h)}\eta(z-ch))| \leq L_gqe^{-\gamma
(t-h)}\eta(z-ch).
$$
Therefore, since $\eta(z) + c\eta'(z) -\eta''(z) >0$ for $z\not=0$ and $\eta(z)
\in (0,1]$, we get
\begin{eqnarray}  \hspace{-2.5cm}
{\mathcal N}w_-(t,z) \leq -q e^{-\gamma t}\left\{\alpha\beta +
\eta(z) + c\eta'(z) -\eta''(z) - \gamma -e^{\gamma h}L_g\right\}< -q e^{-\gamma t}(\gamma +L_g)< 0\nonumber
\end{eqnarray}
for $t\geq 0$. 

\vspace{2mm}

\noindent \underline{Step V}. In view of (\ref{mlem}) and the
monotonicity properties of $g$, we have that
$$
g(w_+(t-h,z-ch))- g(w(t-h,z-ch)) \geq 0, \quad t \in [0,h], \ z
\in \R.
$$
Therefore the difference $\delta(t,z):= w(t,z)- w_+(t,z)$
satisfies the inequalities
\begin{eqnarray}
& &  \hspace{-2.5cm} \delta(0,z) \leq 0,\ |\delta(t,z)| \leq \kappa + q_0^+, \quad
\delta_{zz}(t,z)- \delta_{t}(t,z)-c\delta_{z}(t,z)-\delta (t,z) = \nonumber \\
& &  \hspace{-2.5cm} {\mathcal N}w_+(t,z) - {\mathcal N}w(t,z) + g(w_+(t-h, z-ch))
-
 g(w(t-h, z-ch)) =\nonumber  \\  & &  \hspace{-2.5cm} {\mathcal N}w_+(t,z)+ g(w_+(t-h, z-ch)) -
 g(w(t-h, z-ch))
   \geq 0,  \quad t \in (0,h], \ z \in \R, \ z \not=0; \nonumber  \\
   & &  \hspace{-2.5cm} \frac{\partial \delta(t,0+)}{\partial z}- \frac{\partial \delta(t,0-)}{\partial z}= q\lambda e^{-\gamma t} > 0, \quad t \in (0,h]. \label{di}
\end{eqnarray}
We claim that $\delta(t,z) \leq 0$ for all $t \in [0,h], \ z \in
\R$. Indeed, otherwise there exists $r_0> 0$ such that
$\delta(t,z)$ restricted to any rectangle $\Pi_r= [-r,r]\times
[0,h]$ with $r>r_0$,   reaches its maximal positive value $M >0$
at  at some point $(t',z') \in \Pi_r$.

We claim  that $(t',z')$ belongs to the parabolic boundary
$\partial \Pi_r$ of $\Pi_r$. Indeed, suppose on the contrary, that
$\delta(t,z)$ reaches its maximal positive value at some point
$(t',z')$ of $\Pi_r\setminus \partial \Pi_r$. Then clearly $z'
\not=0$ because of (\ref{di}). Suppose, for instance that $z' >
0$. Then $\delta(t,z)$ considered on the subrectangle $\Pi=
[0,r]\times [0,h]$ reaches its maximal positive value $M$ at the
point $(t',z') \in \Pi \setminus \partial \Pi$.  Then the
classical results \cite[Chapter 3, Theorems 5,7]{PW} shows that
$\delta(t,z) \equiv M >0$ in $\Pi$, a contradiction.

Hence, the usual maximum principle holds for each $\Pi_r, \ r \geq
r_0,$ so that we can appeal to the proof of the
Phragm\`en-Lindel\"of principle from \cite{PW} (see Theorem 10 in
Chapter 3 of this book), in order to conclude that  $\delta(t,z)
\leq 0$ for all  $t \in [0,h], \ z \in \R$.

 But then we can again repeat the above argument on the intervals $[h,2h], \ [2h, 3h], \dots$ establishing that the inequality
 \begin{eqnarray*}\label{eid}
0\leq w(s,z)\leq \phi(z+\epsilon(s))+qe^{-\gamma s}\eta(z), \quad
z\in \R,  \end{eqnarray*}
 actually holds for all $s \geq -h$.  Since $\epsilon(t)$ increases on $\R$, this proves (\ref{mlemm}) with $C=\epsilon(\infty) = \alpha e^{\gamma h}/\gamma$.

Since the same method applied (with  $C=  \alpha e^{\gamma
h}/\gamma$ in (\ref{mlemN})) to the difference $\delta_-(t,z):= w_-(t,z)- w(t,z)$
leads to
\begin{eqnarray*}\label{eid-}  \hspace{-2.5cm}
\phi(z)-qe^{-\gamma s}\eta(z) < \phi(z-\epsilon_1(s))-qe^{-\gamma
s}\eta(z)\leq w(t,z)\leq \kappa, \quad t \geq -h, \ z\in \R,
\end{eqnarray*} the proof of the lemma is completed. \hfill \ter

\begin{remark} \label{R22} It is worthwhile to note that the constants $\gamma, C, q_0^\pm$ 
depend only on the form of $\phi$ in the sense that they will not change if we replace $\phi(z)$  with  a shifted profile $\phi(z+b), \ b \in \R,$ in the statement 
of Lema \ref{lem1}. 
\end{remark}
Due to Remark \ref{R22},   the inequalities (\ref{mlemN}),  (\ref{mlemmN}) can be
presented in the form similar to (\ref{mlem}),  (\ref{mlemm}):
\begin{corollary} \label{corn} Assume that the hypothesis {\rm ({\bf H})} is satisfied.
Then the inequality
\begin{equation*}\label{mlemN2}\phi(z)-q \eta(z) \leq w(s,z)\leq \kappa,  \quad z\in\mathbb{R},\quad s\in [-h,0],
\end{equation*}
with $q \in (0,\sigma]$ implies
 \begin{equation*}\label{mlemmN2}
 \phi(z-Cq)-qe^{-\gamma t} \eta(z) \leq w(t,z)\leq \kappa,
\quad z\in\mathbb{R},\quad t\geq -h.
\end{equation*}
\end{corollary}
{\it Proof. } By Remark \ref{R22}, the statements of Lema \ref{lem1} will not change if we replace $\phi(z)$  with  a shifted profile $\phi(z+b), \ b \in \R$. Taking $b = -Cq$, we complete the proof of Corollary \ref{corn}.  \hfill \ter

\vspace{2mm}

As an immediate consequence of Lemma \ref{lem1} and Corollary
\ref{corn}, we obtain the stability of the wavefront solution
$u(t,x) = \phi(x+ct)$ with respect to the norm $|\cdot|_\lambda$:
\begin{corollary} \label{coco} For every $\epsilon >0$ there exists $\delta >0$ such that $|\phi(\cdot + s_0) - w(s,\cdot)|_\lambda < \delta,$ $s \in [-h,0],$ implies that $|\phi(\cdot  +  s_0) - w(t,\cdot)|_\lambda < \epsilon$ for all $t \geq 0$.
\end{corollary}
{\it Proof. } Without loss of generality, we can assume that $s_0=0$.  From Theorem 1.4 and Proposition 2 from \cite{TPT},   we know that $\phi'(z) =  O(e^{\lambda_{2}z})$ at $-\infty$. This implies that $|\phi'(z)| \leq K\min\{1,\ e^{\lambda_{2}z}\}, z \in \R,$ for some positive $K$. In this way,  for each fixed $p\in \R$,
 $$
0<  \phi'(z+ p) \leq K\min\{1, e^{\lambda_{2}(z+p)}\}\leq Ke^{\lambda_{2}|p|}\min\{1, e^{\lambda_{2}z}\},\  z \in \R.
 $$
Fix $\epsilon >0$ and consider $\delta \in (0, q_0^+)\cap (0,\epsilon/(1+K_1))$, where  $K_1= CKe^{\lambda_{2}Cq_0^+}$.  Next, assume that $|\phi(\cdot) - w(s,\cdot)|_\lambda <
\delta, \ s \in [-h,0]$.  This yields
that $$ \phi(z) - \delta \eta (z) < w(s,z) < \phi(z) +
\delta  \eta (z), \ s \in [-h,0], \ z \in \R,
$$
and therefore, due to Lemma \ref{lem1} and Corollary  \ref{corn},
$$
\phi(z- C\delta) - \delta \eta (z) < w(t,z) < \phi(z+C\delta) + \delta  \eta (z), \ t\geq 0, \ z \in \R.
$$
Now, for some $\hat s \in (0, C\delta)$, it holds
\begin{eqnarray}\nonumber
 \hspace{-.5cm} & &   \phi(z+C\delta) =
\phi(z) + \phi(z+C\delta)- \phi(z) =  \phi(z) + C\delta \phi'(z+ \hat s) \\ 
 \hspace{-.5cm} & &  \leq \phi(z)  + CKe^{\lambda_{2}Cq_0^+}\delta \min\{1,
e^{\lambda_{2}z}\} \leq
  \phi(z)  + K_1\delta \eta(z).\label{it}\end{eqnarray}
After establishing a
similar lower bound  for $\phi(z- C\delta)$, we get
$$
\phi(z) - (K_1+1)\delta  \eta (z) < w(t,z) < \phi(z)
+ (K_1+1)\delta  \eta (z), \  t \geq 0, \ z \in \R,
$$
that is, $|\phi(\cdot) - w(t,\cdot)|_\lambda < \delta
(K_1+1)= \epsilon, \ t \geq 0.$ \hfill \ter

\vspace{2mm}

In addition, Lemma \ref{lem1} yields  the following useful result
\begin{corollary} \label{coro1} Assume that $w_0(s,x)$ satisfies \mbox{\rm ({\bf IC})}. Then
there exist positive $\gamma, \zeta_1$ such that
 \begin{equation}\label{mlemR}  \hspace{-2.5cm}  \phi(z-\zeta_{1})-\sigma e^{-\gamma t}\eta(z)\leq w(t,z)\leq \phi(z+\zeta_1)+q_0^+e^{-\gamma t}
\eta(z+\zeta_1),\ z\in\mathbb{R},\ t\geq -h.
\end{equation}
\end{corollary}
{\it Proof. } \vspace{2mm}
First, we  will show that  inequality (\ref{mlem}) holds for
$w_0(s,z-\zeta_0)$ if we take sufficiently large $\zeta_0$.  Indeed, let $z'$ be such that $\phi(z')+q_0^+\eta(z')=
\kappa$ and define $ \zeta_0$ from
 \begin{eqnarray*}
 Ae^{-\mu \zeta_0} = q_0^+\min\{e^{-\mu z'},e^{(\lambda-\mu)z'}\}.  \end{eqnarray*}
 Then, for all $z \geq z', \ s \in [-h,0]$, it holds that $w(s,z-\zeta_0) \leq 1 \leq \phi(z)+q_0^+\eta(z)$.  Furthermore,  because of the  assumption $(IC2)$ and the inequality $\lambda < \mu$, we have,   for all $z \leq z', \ s \in [-h,0]$,  that $ w(s,z-\zeta_0)\leq Ae^{\mu(z-\zeta_0)} =$
 $$
 q_0^+\min\{e^{-\mu z'},e^{(\lambda-\mu)z'}\}e^{\mu z} \leq q_0^+\min\{e^{\mu (z-z')} ,e^{\lambda z}\} \leq q_0^+ \eta(z) < \phi(z)+ q_0^+ \eta(z).
 $$
Therefore, due to (\ref{mlemm}),
\begin{equation*}
0 \leq w(t,z-\zeta_0)\leq \phi(z+Cq_0^+)+q_0^+e^{-\gamma t} \eta(z),
\quad z\in\mathbb{R},\quad t\geq -h.
\end{equation*}
Hence, setting $\zeta_1= \zeta_0+Cq_0^+$ and using the translation
invariance of equation (\ref{e1}), we obtain the second inequality
in (\ref{mlemR}).

Similarly,  there exists $z''$ such that
\begin{eqnarray} &&
\phi(z+z'')-{q_0^-}\eta(z)\leq 0\leq w(s,z+B)\leq \kappa, \quad
z\leq 0, \ s\in[-h,0]; \nonumber\\ &&
 \phi(z+z'')-{q_0^-}\eta(z)\leq \kappa-\sigma \leq w(s,z+B)\leq \kappa, \quad
z\geq 0, \ s\in[-h,0]. \nonumber
\end{eqnarray}
Hence, by (\ref{mlemmN}) and Remark \ref{R22}, we obtain
\begin{equation*}
 \phi(z+ z''-Cq_0^-)-q_0^-e^{-\gamma t} \eta(z) \leq w(s,z+B)\leq \kappa,
\quad z\in\mathbb{R},\quad t\geq -h.
\end{equation*}
As a consequence, the both inequalities in (\ref{mlemR}) hold if
we take
 $\zeta_1= \zeta_0+C(q_0^++q_0^-)+ |B-z''|$.
\hfill \ter

\begin{remark} \label{dopi} Observe that the hypothesis $(IC3)$ was not used to prove the right-hand side inequality in (\ref{mlemR}). 
\end{remark}
Next, it should be noted that the variable shift $\epsilon(t)$  in
$w_+(t,z)$ was needed only to assure the inequality  ${\mathcal
N}w_+(t,z) \geq 0$ on the finite interval $z-ch+\epsilon(t) \in
[z_1,z_2]$, cf. Step III.   This observation suggests the
following important modification of Lemma \ref{lem1} (where we
will take the same constants $\delta, \gamma>0$ which were  defined in Step II of the
proof of Lemma \ref{lem1}):
\begin{lem}\label{lem2}
Let $w(t,z)$ be a solution of (\ref{wz}), (\ref{wzic}) with  $\tilde w_0(s,z) \in [0,\kappa]$. Take $\delta >0$ as in  (\ref{C1a}) and let $R>ch$ be such that
\begin{eqnarray}  \nonumber 
0\leq w(t,z), \ \phi(z)\leq \delta,  \ \mbox{if} \  z\leq -R+ch,\
t \geq -h, \ \mbox{and}\\ \ |w(t,z)-\kappa|, \
|\phi(z)-\kappa|<\delta\ \ \mbox{if}\  z\geq R-ch, \ t \geq -h.
\nonumber
\end{eqnarray}
Furthermore, suppose that $ w(s,z)\leq \phi(z) + \delta\eta(z)$ \
for all $(s,z) \in [-h,0] \times \R$ and $ w(t,z)\leq \phi(z)$ \
for all $(t,z) \in \R_+ \times [-R-ch,R+ch]$. Then \ $
w(t,z)\leq\phi(z)+ \delta\eta(z)e^{-\gamma t}$\  for all $z\in
\R,\ t\geq 0$.
\end{lem}
{\it Proof. }
Set $\rho(t,z)=w(t,z)-\phi(z)$, then, for some $\xi(t,z)$ lying
between points  $w(t-h,z-ch)$ and $\phi (z-ch)$,
\begin{eqnarray*} \hspace{-2.5cm}
\rho_{t}(t,z)=\rho_{zz}(t,z)-c\rho_{z}(t,z)-\rho(t,z)+g'(\xi(t,z))\rho(t-h,z-ch),
\quad z \in \R,\ t\geq 0.
\end{eqnarray*}
Since $\xi(t,z) \in [0, \delta]$ for $z \leq -R, \ t\geq 0,$ and
$\kappa - \xi(t,z) \in [0, \delta]$ for $z \geq R, \ t\geq 0,$ we find
that $r(t,z):=\delta\eta(z)e^{-\gamma t}$  satisfies
\begin{eqnarray} & & \hspace{-2.5cm}
r_{t}(t,z)-r_{zz}(t,z)+cr_{z}(t,z)+r(t,z)-g'(\xi(z,t))r(t-h,z-ch)=
\nonumber\\ & &  \hspace{-2.5cm} \delta e^{-\gamma t}\left( (1-\gamma) \eta(z) -
\eta''(z) + c \eta'(z) - e^{\gamma h}g'(\xi(t,z))\eta(z-ch)
\right)> 0, \ |z| \geq R, \ t \geq 0. \nonumber
\end{eqnarray}
In addition, by our assumptions,   the piece-wise smooth function
$\delta(t,z): = w(t,z)- (\phi(z)+r(t,z))$ satisfies the
inequalities  $\delta(t,\pm R) \leq 0,\ |\delta(t,z)| \leq 2\kappa +
\delta, \ t \geq 0,$ $z \in\R; $ $ \delta(s,z)\leq 0, \ s \in
[-h,0], \ z \in \R$.  In consequence,
$$\delta_{zz}(t,z)- \delta_{t}(t,z)-c\delta_{z}(t,z)-\delta (t,z) >  -g'(\xi(t,z))\delta(t-h,z-ch) \geq  0,$$
for all $t \in [0,h], \ |z| \geq R$.  By the Phragm\`en-Lindel\"of
principle \cite{PW}, we conclude that   $\delta(t,z) \leq 0$ for
all  $t \in [0,h], \ |z| \geq R$.  Since we also have assumed that
$ w(t,z)\leq \phi(z)$ \ for all $(t,z) \in \R_+ \times
[-R-ch,R+ch]$, we obtain that  $\delta(t,z) \leq 0$ for all  $t
\in [0,h], \ z \in \R$. Finally,  repeating the above arguments on
the intervals $[h,2h], \ [2h, 3h], \dots$, we complete the proof
of the lemma. \hfill \ter

\vspace{2mm}

Finally, before starting with the proof of Theorems \ref{cocoT}  and \ref{main1}, we will establish  the following compactness result.  
\begin{lem} \label{cw}Assume that continuous  function $w: [-h,+\infty) \times \R \to [0,\kappa]$  is a classical solution, for $t>0$, of  equation (\ref{wz}) and  that  $t_j \to +\infty$.
Then there exists a subsequence $\{t_{j_k}\} \subset \{t_j\}$ such
that $w(t_j+s,z)$ converges, uniformly on each rectangle
$[-h,0]\times [-m,m], \ m \in \N$, to  the restriction $w_*(s,z),
\ (s,z) \in [-h,0]\times \R,$ of some entire solution $w_*: \R^2
\to [0,\kappa]$ of   equation (\ref{wz}).
\end{lem}
{\it Proof. } First, we observe that, for each fixed $t >h$, function  $g(w(t-h, z-ch))$ is locally Lipschitz continuous in $z \in \R$ and therefore $w, w_z,w_{zz}$ are H\"older continuous  in $(h,+\infty) \times \R$, cf. \cite[Theorem 1]{Kn}.  Next, fix an arbitrary positive $T>2h+2$ and $m\in \N$ and consider, for $t_j >T+2h$,  solutions
$w_j(t,z) = w(t_j+t,z), \ (t,z) \in D_+:= [-T,T]\times
[-m-1,m+1+ch],$ of the equation
\begin{eqnarray*}
w_{t}(t,z)= w_{zz}(t,z)-cw_{z}-w(t,z)+g_j(t,z),
\end{eqnarray*}
where $g_j(t,z):= g(w_j(t-h, z-ch))$. We claim that, for each
$\alpha \in (0,1)$,  there exists a positive $K$ depending only on
$m, T, \alpha$ such that the H\"older norms
$$
|g_j|_\alpha^D=  \sup_{(t,z)\in D}|g_j(t,z)|+
\sup_{(t,z)\not=(s,x)\in
D}\frac{|g_j(t,z)-g_j(s,y)|}{(|x-z|^2+|t-s|)^{\alpha/2}}
$$
are uniformly bounded  in $D:=[-T+1+h,T]\times [-m,m]$ by $K$
(i.e. $|g_j|_\alpha^D \leq K$ for all $j$. Observe  that
$|g_j|_\alpha^{D_+}$ is finite due to \cite[Theorem 1]{Kn}).  In
fact,  since $g$ satisfies the Lipschitz condition on
$[0,\kappa]$,  it suffices to establish the uniform boundedness of
$|w_{j}|_\alpha^{D_1}$ in a bigger domain $D_1:=[-T+1,T]\times
[-m,m+ch]\subset D_+$. Obviously,  $w_j$ solves in $D_+$ the
initial-boundary value problem $w= w_j|_{\partial D_+}$ where
$w_j|_{\partial D_+}$ denotes the restriction of $w_j$ on the
parabolic boundary $\partial D_+:=  \{-T\} \times [-m-1,m+1+ch]
\cup [-T,T]\times \{-m-1, m+1+ch\}$ of $D_+$.  Let $\rho:
[-T,T]\to [0,1]$ be some nondecreasing smooth function such that
$\rho([-T,-T+0.25])=0,$  $\rho([-T+0.75,T])=1$.  Then $w_j=
w_{j,1}+ w_{j,2}$ where $w_{j,1}$ is the solution of the
initial-boundary value problem $w= 0|_{\partial D_+}$ for the
equation
$$
w_{t}(t,z)= w_{zz}(t,z)-cw_{z}-w(t,z)+\rho(t)g_j(t,z),
$$
and  $w_{j,2}$ solves the initial-boundary value problem $w=
w_j|_{\partial D_+}$ for the equation
$$
w_{t}(t,z)= w_{zz}(t,z)-cw_{z}-w(t,z)+(1-\rho(t))g_j(t,z).
$$
Next, since $|g_j(t,z)|\leq \kappa$ for all $(t,z) \in D_+, j \in
\N,$ a priori estimate (of the type $1+\delta$) established in
\cite[Theorem 4, Chapter 7]{AF} guarantees that
$|w_{j,1}|_\alpha^{D_+} \leq K_1, j \in \N$, where $K_1$ depends
only on $m, T, \alpha$. As consequence, since $\sup_{D_+}|w_j|, \ j \in \N,$
are uniformly bounded by $\kappa$, we deduce that
$\sup_{D_+}|w_{j,2}|=\sup_{D_+}|w_j-w_{j,1}|, \ j \in \N,$  are  also uniformy
bounded. In addition,  $(1-\rho(t))g_j(t,z) =0$ in $
[-T+0.75,T]\times [-m-1,m+1+ch],$  so that  we can invoke the
interior Schauder estimates (see, e.g,  \cite[Theorem 5, Chapter
3]{AF}) in order to deduce  that  $|w_{j,1}|_\alpha^{D_1} \leq
K_2, j \in \N,$ where  $K_2>0$ depends only on  $\alpha$ and
$K_1$. Hence, $|w_{j}|_\alpha^{D_1} \leq K_1+K_2, j \in \N,$ and
therefore  $|g_j|_\alpha^D \leq K:= L_g(K_1+K_2)$ for all $j$.

Applying now Theorem 15 from \cite[Chapter 3]{AF},  we conclude
that there exists a subsequence $\{t_{j_k}\} \subset  \{t_j\}$
such that $w_{j_k}(t,z)$ converges, uniformly on
$[-T+2+h,T-1]\times [-m+1,m-1],$ to the classical solution
$w_{T,m}: [-T+2+2h, T]\times [-m+1, m-1] \to [0,\kappa]$  of
equation (\ref{wz}).   Finally, considering $m,T \to +\infty$ and
applying a standard diagonal argument, we can assume that
$w_{j_k}(t,z)$ converges,  uniformly on compact subsets of $\R^2$
to an entire classical solution $w_*: \R^2\to [0,\kappa]$  of  the
functional differential equation  (\ref{wz}). Observe that the
arguments used to estimate $|w_{j}|_\alpha^{D_1}$ can be also
applied without changes to $w_*$ so that $|w_{*}|_\alpha^{D_1} \leq
K_1+K_2$ with the same $K_1, K_2$.
\hfill \ter

 \begin{remark}\label{ols}
 Due to Lemma \ref{cw}, we can define $\omega$-limit set $\omega(w_0)$ which consists from the restrictions $w_*(s,z),$ $(s,z) \in [-h,0]\times \R,$  of all possible entire limit solutions $w_*=\lim_{k \to +\infty} w_{j_k}$ to  (\ref{wz}) (which are obtained by considering all possible sequences $\{t_j\}$ converging to  $+\infty$ in this lemma). Since each $w_*$ is an entire solution, the set $\omega(w_0)$ is invariant. Furthermore,  since $|w_{*}|_\alpha^{D} \leq K_1+K_2$ where  $K_1, K_2$ depend only on $D$ and $\alpha$, the set $\omega(w_0)$ is pre-compact with respect to the topology of the uniform convergence on bounded subsets
 of $[-h,0]\times \R$.  Actually  $\omega(w_0)$ is compact in the mentioned topology since each element of
 $\omega(w_0)$ can uniformly (on bounded sets) approximated by  $w_j$.
 \end{remark}
\begin{thm} \label{cocos}Assume that $u = \phi(x +c_*t), \ c_* > c_\#$,  is a pushed traveling front  to equation (\ref{e1}). If initial  function $w_0$ satisfies all conditions {\rm ({\rm \bf IC})} then, for some $z_0\in \R$,   the classical solution $w= w(t,z)$ of the initial value problem (\ref{wz}), (\ref{wzic}) asymptotically converges to  a shifted front profile:
\begin{equation}\label{conv}
\lim_{t \to \infty }|w(t,\cdot ) -   \phi(\cdot + z_0)|_{\lambda}
=0.
\end{equation}
\end{thm}
In order to prove the above theorem, instead of looking for an
appropriate Lyapunov functional (as it was done in \cite{FML,
roth}) for  functional differential equation
(\ref{wz}), we will use the Berestycki and Nirenberg method of the
sliding solutions as well as some ideas of the   approach
developed by Ogiwara and Matano in \cite{OM}.

{\it Proof. } By Corollary \ref{coro1}, Lemma \ref{cw} and Remark \ref{ols}, solution $w= w(t,z)$ of the initial value problem (\ref{wz}), (\ref{wzic}) has a compact invariant $\omega-$limit set $\omega(w_0)$ such that for some fixed $\zeta_1$, it holds
\begin{equation}\label{ps}
\phi(z-\zeta_{1})\leq w_*(0,z)\leq \phi(z+\zeta_1), \quad z\in\R,
\ \mbox{ for each} \  w_* \in  \omega(w_0).
\end{equation}
 Then the set
$$
A = \{a\in \R: w_*(0,z)\leq \phi(z+a), \ z \in \R, \ \mbox{for
each } \ w_* \in  \omega(w_0)\}
$$
contains $\zeta_1$ and has $-\zeta_1$ as its lower bound.
Therefore $\hat a= \inf A$ is a well defined finite number.   Due
to continuity of $\phi$, we have that  $\hat a \in A$ so that
$$
w_*(0,z)\leq \phi(z+\hat a), \ z \in \R, \ \mbox{for each } \ w_*
\in  \omega(w_0).
$$
In fact, since $\omega(w_0)$ is an invariant set, we have that
$w_*(t,z) \leq \phi(z+\hat a), \ z \in \R, \ t \in \R$.  Suppose
now for a moment that $w_*(0,z')= \phi(z'+\hat a)$ for some finite
$z'$ and some $w_* \in \omega(w_0)$.  Therefore, since $g$ is an increasing function, the strong
maximum principle yields $w_*(t,z) \equiv \phi(z+\hat a)$ for all
$t\leq 0, \ z \in \R$. In particular, $w_*(0,z) \equiv \phi(z+\hat
a)$ so that, for some sequence $t_n\to +\infty$, it holds that
$w(t_n+s,z) \to \phi(z+\hat a)$  uniformly with respect to $s \in [-h,0]$ and $z$ from compact subsets of
$\R$.  In addition,  Corollary \ref{coro1} allows to evaluate the difference  $|w(t_n+s,z) - \phi(z+\hat a)|/\eta(z)$ in 
some fixed  neighbourhood of the endpoints $z=-\infty$ and $z= +\infty$ and to conclude that $w(t_n+s,z) \to \phi(z+\hat a), n \to +\infty$, in the norm $|\cdot
|_\lambda$ and uniformly with respect to $s \in [-h,0]$. By Corollary \ref{coco}, the latter convergence
implies  (\ref{conv}) with $z_0 =\hat a$ that completes  the proof of the theorem in
the case when $w_*(0,z')= \phi(z'+\hat a)$ holds for some finite
$z'$.

In this way,  we are left to consider the situation when
\begin{equation}\label{psi}
w_*(0,z)< \phi(z+\hat a), \ z \in \R, \ \mbox{for each } \ w_* \in
\omega(w_0).
\end{equation}
In virtue of (\ref{ps}), for any given $\delta >0$, we can find $R>3ch+1$ sufficiently large
to have, for all $w_* \in  \omega(w_0)$,
$$
w_*(0,z) <  \phi(z+\hat a) < \delta, \ \mbox{for} \ z \leq
-R+ch+1, \quad \phi(z+\hat a) > w_*(0,z) > \kappa- \delta, \
\mbox{for} \ z \geq R -ch-1.
$$
Then, using (\ref{psi}) and the compactness of the set
$$\{w_*(0,\cdot): [-R+ch+1, R-ch-1] \to [0,\kappa] , \quad  w_*
\in \omega(w_0)\} \subset C[-R+ch+1, R-ch-1], $$ we deduce the
existence of $\varsigma\in (0,1)$ such that
\begin{equation*}\label{psi2}\hspace{-2.5cm}
w_*(0,z)< \phi(z+\hat a-\varsigma), \ z \in [-R+ch+1, R-ch-1], \
\mbox{for each } \ w_* \in  \omega(w_0).
\end{equation*} It is clear that
$$
\phi(z+\hat a) < \kappa < \phi(z+\hat a-\varsigma) + \delta, \quad z
\geq R-ch.
$$
Without the loss of generality, we also can suppose that
$\varsigma\in (0,1)$ is such that
$$
\phi(z+\hat a) < \phi(z+\hat a-\varsigma) + \delta e^{\lambda z},
\quad z \leq -R+ch+1.
$$
Indeed, observe that $\phi'(z) \leq Ce^{\lambda_2z}, \ z \leq 0,$
and therefore, for some $\xi \in (z+\hat a-\varsigma,z+\hat a),$
$$
\phi(z+\hat a) - \phi(z+\hat a-\varsigma) = \phi'(\xi)\varsigma \leq
Ce^{\lambda_2 (z+\hat a)} \varsigma \leq \delta e^{\lambda z},\ z \leq 0,
$$
once $\varsigma \leq e^{-\lambda_2 \hat a} \delta/C $. Hence,
invoking again the invariance property of $\omega(w_0)$, we can
conclude that for each $w_* \in \omega(w_0)$ it holds
$$
w_*(t,z) \leq \phi(z+\hat a-\varsigma) + \delta \eta(z), \quad t \in
\R, \ z \in \R, $$ and
$$ w_*(t,z) \leq \phi(z+\hat a-\varsigma), \ z \in [-R+ch+1, R-ch-1], \ t \in \R.
$$
By Lemma \ref{lem2}, this yields $ w_*(t,z)\leq\phi(z+\hat
a-\varsigma)+ \delta\eta(z)e^{-\gamma t}, \ t \geq 0, \ z \in \R,$
where $\hat a, \varsigma, \gamma$ do not depend on the particular
choice of $w_* \in \omega(w_0)$.  In consequence,  since
$w_*(t,z)$ is an entire solution,  we obtain that actually
$w_*(0,z) \leq \phi(z+\hat a-\varsigma), \ z \in \R,$ for all
$w_*\in \omega(w_0)$. This contradicts to the definition of $\hat
a$ and shows that the case (\ref{psi}) 
can not happen.\hfill \ter

\section{Proof of Proposition \ref{mainP2}} \label{sub2}
First, observe that for each $g$ satisfying the assumptions of  Proposition \ref{mainP2}, we can find a {\it monotone} function 
$g_1:[0,\kappa]\to [0,\kappa]$ possessing all the properties of $g$ and such that $g_1(x)\leq g(x)$. Therefore, in view of the comparison principle, it will not restrict the generality if  we will assume additionally  the monotonicity of  $g$.  

\noindent Here, we   follow an approach, proposed by Aronson and Weinberger in \cite[Theorem 3.1]{AW},  and based on  the maximum principle.  In the mentioned work, it was established, 
for every $\epsilon \in (0,\kappa)$ and appropriate $b_\epsilon >0$, the existence of a positive solution $q = q(x) \leq \epsilon$ to  the Dirichlet boundary value problem
\begin{eqnarray*}\label{eq2o}
& & q''(x)-q(x)+g(q(x))=0, \  x\in I_{\epsilon}:=(0,b_{\epsilon}),\\
& & q(0)=q(b_{\epsilon})= 0.
\end{eqnarray*}
Since we are interested in the asymptotic behavior of $u(t,x)\geq 0$ and $w_0(s,x) \not\equiv0$,   without loss of
generality, 
due to the  strong maximum principle we can suppose that $w_{0}(s,x)>0$ for all 
$(s,x)\in[-h,0]\times\R$. But then we can choose $\epsilon >0$ small enough to have  $q(x)\leq u(x,s)$ for all $x \in {I}_{\epsilon}$, $ s \in [-h,0]$.  Let $\chi_A$ denote the characteristic function of  subset $A \subset \R$. Consider  solution $u =  u_\epsilon (t,x)$ of the initial value problem  $u_\epsilon (s,x) = \chi_{I_\epsilon}(x)q(x),$ $s \in [-h,0],\ x \in \R$, to equation (\ref{e1}).  Then the difference  
$\delta(t,x)=q(x)- u_\epsilon (t,x)$ satisfies $\delta(t,0) \leq 0, \ \delta(t,b_{\epsilon}) \leq 0, \ t \geq 0,$ and 
$$
\delta_{t}(t,x)-\delta_{xx}(t,x)+\delta(t,x)=g(q(x))- g(u_\epsilon (t-h,x))\leq
0, \quad (t,x)\in[0,h]\times I_{\epsilon}.
$$
Hence, by 
the maximum principle, $u_{\epsilon}(t,x)\geq q(x)$ on
$[0,h]\times I_{\epsilon}$. Repeating the same argument on $[h,2h]\times I_{\epsilon}$, we
obtain that $u_{\epsilon}(t,x)\geq q(x)$ for all $(t,x) \in [h,2h]\times  I_{\epsilon}$.  It is clear that 
this procedure yields the inequality   
$q(x) \leq u_{\epsilon}(t,x) < 1$ in $[0,+\infty)\times I_{\epsilon}$.  
But then, since for each positive $l$, it holds that $\chi_{I_\epsilon}(x)q(x) = u_{\epsilon}(s,x) \leq u_{\epsilon}(s+l,x),$  $s \in [-h,0], \ x \in \R,$ we can 
use the 
Phragm\`en-Lindel\"of principle, in order to 
conclude that  $u_{\epsilon}(t+l,x)\geq u_{\epsilon}(t,x)$ for all $(t,x) \in [0,h]\times \R$. Similarly to the above analysis, step by step, we can extend the latter inequality  for all $(t,x) \in \R_+\times \R$. Thus, for each fixed $x \in \R$, 
$u_{\epsilon}(t,x)$ is a non-decreasing bounded function of $t \geq 0$.  Let $u_{\epsilon}(x) = \lim_{t \to +\infty} u_{\epsilon}(t,x)$, then $u_{\epsilon}(x) \in (0, \kappa]$ for every $x \in \R$.  

Now, a direct application of Lemma \ref{cw} shows that $u_{\epsilon}(x)$ solves 
$$
u''(x)-u(x)+g(u(x))=0, \quad x\in \R
$$
while the convergence $u_{\epsilon}(x) = \lim_{t \to +\infty} u_{\epsilon}(t,x)$ is uniform on compact subsets of $\R$. 
Since $g(u)-u >0$ on $(0,\kappa)$, the function $u_{\epsilon}(x)$ cannot take (local) minimal values in $(0, \kappa)$. This  implies the existence of $u_{\epsilon}(\pm\infty)\in \{0,\kappa\}$.  In other words, $u(x)$ is a positive stationary traveling wave solution 
of equation (\ref{e1}) considered with $h=0$.  It is well known \cite{SEDY} that this is possible only when $u_{\epsilon}(x) \equiv \kappa$. 

Finally, we complete the proof by observing that,  due to  the maximum principle, it holds  
$u_{\epsilon}(t,x) \leq u(t,x)$ on $[0,\infty)\times\R$.

\section{Proof of Theorem \ref{main2}: auxiliary results}\label{sub3}
In  Sections \ref {sub3} and \ref {sub4}, we are always assuming that all the conditions of Theorem \ref{main2} are
satisfied (recall also that, by simplifying the notation, we write $c$ instead of $c_*$). The proof of this theorem will follow from a series of lemmas.  In the first of them 
we improve the asymptotic relation $u(t,0) =\kappa + o(1)$ at $+\infty$ known from Proposition \ref{mainP2}.   As we show below,  this convergence is actually of the exponential type. 
 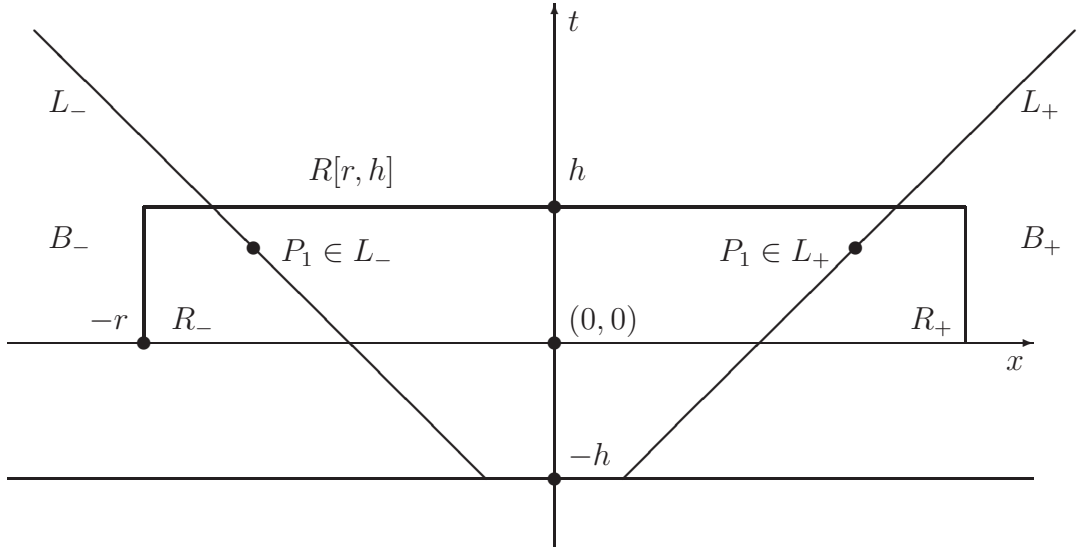
\begin{figure}[h]\label{FFF2}
 \setlength{\unitlength}{1.8cm}
\begin{picture}(1,4)(-0.5,-0.5)
\put(-0.5,1){\vector(1,0){7.5}}
\put(6.8,0.8){$x$}
\put(3.6,3.3){$t$}
\put(0.1, 1.1){$-r$}
\put(-0.2, 1.7){$B_-$}
\put(6.9, 1.7){$B_+$}
\put(-0.2, 2.7){$L_-$}
\put(6.9, 2.7){$L_+$}
\put(0.7, 1.1){$R_-$}
\put(1.7, 2.2){$R[r,h]$}
\put(6.1, 1.1){$R_+$}
\put(0.5,1){\circle*{0.1}}
\put(3.5,1){\circle*{0.1}}
\put(3.5,2){\circle*{0.1}}
\put(3.5,-0){\circle*{0.1}}
\put(5.7,1.7){\circle*{0.1}}
\put(1.3,1.7){\circle*{0.1}}
\put(4.7,1.6){$P_1\in L_+$}
\put(1.5,1.6){$P_1\in L_-$}
\put(3.6,0.1){$-h$}
\put(3.6,1.1){$(0,0)$}
\put(3.6,2.2){$h$}
\put(3.5,-0.5){\vector(0,1){4}}
\put(-0.5,0){\line(1,0){7.5}}
\thicklines
\put(0.5,1){\line(0,1){1}}
\put(6.5,1){\line(0,1){1}}
\put(0.5,2){\line(1,0){6}}
\put(4,0){\line(1,1){3.3}}
\put(3,0){\line(-1,1){3.3}}
\end{picture}
\caption{\hspace{0.1cm} Domains $R[r,h], R_1, R_\pm \subset B_\pm$ and lines $L_\pm$.}
\end{figure}
\begin{lem}\label{x0} Assume 
that 
${\lambda_1}<-{\lambda}_{3}$  where 
$\lambda_{3}$ stands for a unique negative zero of the characteristic function
$
\chi_\kappa(z,c):= z^{2}-cz-1+g'(\kappa)e^{-zch}
$.   If ${u}(t,x)$ solves (\ref{e1}), (\ref{e2}) with $w_0 \not\equiv 0$, 
then 
there exist numbers  $q,\nu>0$ such that
\begin{eqnarray}\label{eq14} u(t,0) \geq 
\kappa-qe^{-\nu t}\quad \mbox{for all}\ \ t\geq 0. 
\end{eqnarray}
\end{lem}
{\it Proof. } 
First, we fix a positive $\lambda \in (\lambda_1, -\lambda_3)\cap(\lambda_1, \lambda_2) $ and $\gamma < \min \{ c\lambda, \gamma_1^*\}$ such 
that 
\begin{eqnarray*}\label{C1aa} \hspace{-.5cm} 
-\lambda^2+c\lambda +1-\gamma-g'(\bar s)e^{-\lambda ch+\gamma
h}>0 \ \mbox{for} \ \bar s < \delta, 
\end{eqnarray*}
where $\delta < \delta_1^*, \gamma_1^*, z_1 <0<z_2$ are defined in Steps I, II of Lemma \ref{lem1}. 
Following \cite{FML},  we will construct a  sub-solution to (\ref{e1}) of the form 
$$
{u}_{-}(t,x)={\phi}_{+}(t,x)+{\phi}_{-}(t,x)-\kappa-{q}(t,x),
$$
where ${\phi}_{\pm}(t,x)={\phi}(\pm
x+{c}t-\epsilon(t)), \ {q}(t,x)=\gamma e^{-\gamma
t}\theta(t,x)$  with $\theta(t,x) \leq 1$ are defined by  
\begin{eqnarray*}  \hspace{-2.5cm}
\theta(t,x):=\eta(-|x|+{c}t-\epsilon(\infty)-z_{1})=\left\{
\begin{array}{ll} \hspace{-2mm}
     &  \hspace{-5mm} e^{ {\lambda}(-x+ {c}t-\epsilon(\infty)-z_{1})},  \  \mbox{if} \  {(t,x)\in B_{+}}, \\
     &  \hspace{-5mm}  e^{ {\lambda} (x+ {c}t-\epsilon(\infty)-z_{1})}, \  \mbox{if} \ 
     {(t,x)\in B_{-},}\\
   &  \hspace{-5mm}  1, \  \mbox{if} \  {(t,x)\in[-h,\infty)\times\mathbb{R}\setminus(B_{+}\cup B_{-})}, 
\end{array}%
\right.
\end{eqnarray*}
 $$B_{\pm}:=[-h,\infty)\times\mathbb{R}\cap\{(t,x):\mp x+ {c}t-\epsilon(\infty)<z_{1}\},$$ 
 $$L_{\pm}:=[-h,\infty)\times\mathbb{R}\cap\{(t,x):\mp x+ {c}t-\epsilon(\infty)
 =z_{1}\}, $$
with an
appropriate $\epsilon(t)$ satisfying $\epsilon'(t)>0,$ $
\epsilon(t)<0$. Then  {$\epsilon(\infty)+z_1 < -ch$ and therefore $B_+\cap B_-=\emptyset$}. See also Figure 1.  Set 
$$
\mathcal{N}_1 {u}_{-}(t,x):=
 ( {u}_{-})_{t}(t,x)-( {u}_{-})_{xx}(t,x)+ {u}_{-}(t,x)- {g}(u_-(t-h,x)),$$
 $$
\tilde \phi_\pm(t,x):= {\phi}(\pm x+ {c}(t-h)-\epsilon(t)) <  \phi_\pm(t-h,x)
.$$ 
Since ${u}_{-}(t,x) = {u}_{-}(t,-x)$, it holds that $\mathcal{N}_1 {u}_{-}(t,x)= \mathcal{N}_1 {u}_{-}(t,-x)$. In view of monotonicity of
$ {g}$ and $ {\phi}$, we have
\begin{eqnarray}
& &\hspace{-25mm} \mathcal{N}_1 {u}_{-}(t,x)\leq  {g}(\tilde \phi_+(t,x))+ {g}(\tilde \phi_-(t,x)) - {g}\left(\tilde \phi_+(t,x)+ \tilde \phi_-(t,x)-\kappa- {q}(t-h,x)\right)\nonumber\\
&&\hspace{-25mm} -\epsilon'(t)[ {\phi}'(x+ {c}t
-\epsilon(t))+ {\phi}'(-x+ {c}t
-\epsilon(t))]-\kappa- {q}(t,x)+ {q}_{xx}(t,x)- {q}_{t}(t,x).\nonumber
\end{eqnarray}

\noindent \underline{Claim I:  $\mathcal{N}_1 {u}_{-}(t,x) = \mathcal{N}_1 {u}_{-}(t,-x) <0$ for $x \geq 0, \ t >0,\ (t,x) \not\in L_\pm$}. 

\noindent   By  Step II
of Lemma \ref{lem1}, 
 $-x+ {c}(t-h)-\epsilon(t) \geq z_2$ implies   $\kappa-\delta/2<\tilde \phi_{-}(t,x)$.
 Since $x\geq 0$, we also have $\kappa-\delta/2<\tilde \phi_{-}(t,x)\leq \tilde \phi_{+}(t,x)$.
  By Step I of Lemma \ref{lem1}, for $\gamma \in (0,\gamma_1^*), \ \gamma < \sigma-\delta/2$, 
\begin{eqnarray*}\label{eq14.1}\hspace{-0.5cm}
  g(\tilde \phi_{-}(t,x))-  g(\tilde
\phi_{-}(t,x)-[\kappa-\tilde \phi_{+}(t,x)+  q(t-h,x)])\nonumber \\ \leq
e^{-\gamma h}(1-2\gamma)[\kappa-\tilde \phi_{+}(x,t)+ 
q(t-h,x)].
\end{eqnarray*}
Hence, since  $\theta(t,x)$ is non-decreasing in
$t$, we have, {for $t >0$, } that $ \mathcal{N}_1u_{-}(t,x)\leq$
\begin{eqnarray*}&&\hspace{-25mm}
(1-2\gamma)e^{-\gamma
h}[(\kappa-\tilde \phi_{+}(t,x))+  q(t-h,x)]+  g(\tilde \phi_{+}(t,x))-\kappa-  q(t,x)+  q_{xx}(t,x)-  q_{t}(t,x)  \\
&&  \hspace{-2.5cm}\leq -\tilde \phi_{+}(t,x)+(1-2\gamma)\gamma e^{-\gamma
t}\theta(t-h,x)+  g(\tilde \phi_{+}(t,x)) -  q(t,x)+  q_{xx}(t,x)-  q_{t}(t,x)\\
&& \hspace{-25mm}   \leq g(\tilde \phi_{+}(t,x))-\tilde \phi_{+}(t,x)
+  q(t,x)  \left\{
\begin{array}{ll}  &  \hspace{-5mm}
     \lambda^2 - {c}  \lambda  -1 +\gamma + (1-2\gamma){e^{ -  \lambda  {c}h}}, \ \mbox{if} \ {(t,x)\in B_{+}},\\  &  \hspace{-5mm}
     -2\gamma, \ \mbox{if} \
{(t,x)\in[0,\infty)\times\mathbb{R}^{+}\setminus B_{+}}.
\end{array}%
\right.
\end{eqnarray*}
On the other hand, it is known (see e.g. \cite[Remark 1]{TPT}) that, for  some $C>0$, it holds 
\begin{eqnarray}\label{eq14.2}
0\leq \kappa-\tilde \phi_{+}(t,x)\leq
Ce^{- {\lambda}_{3}\epsilon(t)}e^{ {\lambda}_{3}(x+ {c}t)}, 
\qquad t \geq -h, \ x \in \R.
\end{eqnarray}
This implies that, for $t >0, x \geq 0$,  $-x+ {c}(t-h)-\epsilon(t) \geq z_2$,  it holds that
$$\mathcal{N}_1u_{-}(t,x)\leq Ce^{ -{\lambda}_{3}\epsilon(\infty)}e^{ {\lambda}_{3}x}e^{ {\lambda}_{3} {c}t}$$
\begin{eqnarray*}&& \hspace{-25mm}
+ \gamma e^{-\gamma t} \left\{
\begin{array}{ll} &  \hspace{-5mm}
    e^{  \lambda(-x+ {c}t-\epsilon(\infty)-z_{1})}[  \lambda^2 - {c}  \lambda  -1 +\gamma + (1-2\gamma)e^{-  \lambda  {c}h}],  \ \mbox{if} \  {(t,x)\in B_{+}},\\ &  \hspace{-5mm}
     -2\gamma, \ \mbox{if} \  
{(t,x)\in[0,\infty)\times\mathbb{R}^{+}\setminus B_{+}},
\end{array}%
\right. 
\end{eqnarray*}
\begin{eqnarray*} \hspace{-25mm}
\leq e^{-\gamma t} \left\{
\begin{array}{ll} &  \hspace{-5mm}
   e^{-\lambda x}[\gamma(\lambda^2 - {c}  \lambda  -1 +\gamma + (1-2\gamma)e^{-  \lambda  {c}h})+Ce^{ -{\lambda}_{3}\epsilon(\infty)}],  \ \mbox{if} \  {(t,x)\in B_{+}},\\ &  \hspace{-5mm}
     -2\gamma^2 +Ce^{ -{\lambda}_{3}\epsilon(\infty)},  \ \mbox{if} \ 
{(t,x)\in[0,\infty)\times\mathbb{R}^{+}\setminus B_{+}}.
\end{array}%
\right.
\end{eqnarray*}
As a consequence, there
exists  large negative $\epsilon(\infty)$
(depending on $\gamma$ and $\lambda_3$) such that 
$$
\mathcal{N}_1  u_{-}(t,x)
<0\ \mbox{for} \
 t >0, \ -x+ {c}(t-h)-\epsilon(t) \geq z_2, \ (t,x) \not\in L_+. 
$$
Next, if   $-x+ {c}(t-h) - \epsilon(t) \leq z_1$ then $0\leq\tilde
\phi_{-}(t,x)\leq \delta/2$ and $(t,x) \in B_+$. Thus $\theta(t-h,x)=e^{ 
\lambda(-x+ {c}t- {c}h-\epsilon(\infty)-z_{1})}$ and, for
some $\bar s < \delta/2$,
\begin{eqnarray*}\hspace{-25mm}\nonumber
  g(\tilde\phi_{-}(t,x))- 
g(\tilde\phi_{-}(t,x)-[\kappa-\tilde\phi_{+}(t,x)+ 
q(t-h,x)]) =   g'(\bar{s})[\kappa-\tilde\phi_{+}(t,x)+ 
q(t-h,x)].
\end{eqnarray*}
Thus, recalling that $z_1<0$, for large $\epsilon(\infty)<0$ (which depends
on $\gamma$ and $\lambda_3$),  we get  \begin{eqnarray*}
&& \hspace{-25mm}\mathcal{N}_1  u_{-}(t,x)\leq
  g'(\bar{s})[\kappa-\tilde \phi_{+}(t,x)+  q(t-h,x)]-  q(t,x)+  q_{xx}(t,x)-  q_{t}(t,x)\nonumber\\
 && \hspace{-25mm} \leq 
      \gamma e^{-\gamma t}e^{  \lambda(-x+ {c}t-\epsilon(\infty)-z_{1})}[  \lambda^{2}- {c}  \lambda-1+\gamma+  g'(\bar{s})e^{(\gamma- {c}  \lambda)h}]+    g'(\bar{s})[\kappa-\tilde \phi_{+}(t,x)]\nonumber\\
     &&  \hspace{-25mm} \leq  \gamma e^{-\gamma t}e^{-  \lambda x}[  \lambda^{2}- {c}  \lambda-1+\gamma+  g'(\bar{s})e^{(\gamma- {c}  \lambda)
     h}]+  g'(\bar{s})Ce^{ {\lambda}_{3}(x+ ct-\epsilon(t))} \nonumber\\
     &&  \hspace{-25mm}  \leq  e^{-  \lambda x}e^{-\gamma t}\left\{\gamma[\lambda^{2}- {c}  \lambda-1+\gamma+  g'(\bar{s})e^{(\gamma- {c}  \lambda)
     h}]+  g'(\bar{s})Ce^{- {\lambda}_{3}\epsilon(\infty)}\right\}<0.
\end{eqnarray*}
Finally, consider  $z_1 \leq  -x+ {c}(t-h) - \epsilon(t) \leq z_2$. Recall that  
$\beta>0$ defined in Step II of Lemma \ref{lem1} depends only on $\delta, \phi$ and satisfies $\beta < \min_{\zeta\in[z_{1},z_{2}+ch]}\phi'(\zeta)$. Therefore, if we take 
$\epsilon'(t)= \alpha \gamma e^{-\gamma t}$ for some $\alpha>0$, then 
$$
|  g(\tilde\phi_{-}(t,x))- 
g(\tilde\phi_{+}(t,x)+\tilde\phi_{-}(t,x)-\kappa- 
q(t-h,x))|\leq L_{ 
g}[Ce^{- {\lambda}_{3}\epsilon(t)}e^{ {\lambda}_{3}(x+ {c}t)}+ 
q(t-h,x)].
$$
In consequence, if $\alpha$ is sufficiently large then  
$\mathcal{N}_1 
u_{-}(x,t)\leq CL_{  g}
e^{ {\lambda}_{3}(-\epsilon(t)+x+ {c}t)}$
\begin{eqnarray*} \hspace{-25mm}
+\left\{
\begin{array}{lll} &  \hspace{-5mm}
      \gamma e^{-\gamma t}\{-\alpha\beta+e^{  \lambda(-x+ {c}t-\epsilon(\infty)-z_{1})}[  \lambda^{2}- {c}  \lambda-1+\gamma+L_{  g}e^{(\gamma-  \lambda  {c})h}]\},   \ \mbox{if} \  {(t,x)\in B_{+}}, \\ &  \hspace{-5mm}
 \gamma  e^{-\gamma t}[-\alpha\beta+\gamma -1 +L_{g}e^{\gamma h}],   \ \mbox{if} \ 
{(t,x)\in[0,\infty)\times\mathbb{R}^{+}\setminus B_{+}},
\end{array}%
\right. \leq
\end{eqnarray*}
$$
e^{-\gamma t}\left\{ \gamma [-\alpha\beta+L_{g}e^{\gamma h}]+CL_{g}e^{- {\lambda}_{3}\epsilon(\infty)}
e^{\lambda_{3}x}\right\}<0,   \  \mbox{for} \  (t,x)\in \mathbb{R}^{+}\times\mathbb{R}^{+}.
$$
\noindent \underline{Claim II:  There exists $t_0>0$ such that $u_{-}(s,x)\leq u(s+t_0,x)$ for $x \in \R, \ s \in [-h,0]$}. 
 
\noindent 
Since $\lambda_{2}>\lambda$, 
there exists $r_0>0$ depending  on $\epsilon(-h), \epsilon(\infty), z_1$ such that, for $s\in[-h,0]$, 
\begin{eqnarray*}  \hspace{-15mm}
   u_{-}(s,x)&\leq&   {\phi}(-|x|+  {c}s-\epsilon(s))- \gamma \eta(-|x|+  {c}s-z_{1}-\epsilon(\infty)) <0 \ \mbox{if} \
 |x|\geq r_0. 
\end{eqnarray*}
Clearly,   $u_{-}(s,x) < \kappa$  for all $|x|\leq r_0, s\in[-h,0]$ and therefore, by Proposition \ref{mainP2}, 
$u_{-}(s,x) < u(t_{0}+s,x)$, $|x|\leq r_0, s\in[-h,0]$, for an appropriate $t_0>0$. 

\vspace{2mm}

Claims I and II allow to complete the proof of Lemma \ref{x0}.  First, for $r>0$, consider rectangle
$R[r,h]=[0,h]\times[-r,r]$. Set $\delta(t,x):=   {u}_{-}(t,x)-  {u}(t+t_0,x)$, the   function $\delta(t,x)$ is smooth in $[-h,+\infty)\times \R\setminus\{L_-\cup L_+\}$ (in particular, in  the regions $R_\pm = R[r,h]\cap B_\pm$, $R_1 = R[r,h]\setminus(\bar R_+ \cup \bar R_-)$). Since $\delta(s,x) \leq 0$ in $[-h,0]\times \R$ and
$$
\delta_{t}(t,x)-\delta_{xx}(t,x)+\delta(t,x)\leq g(u_-(t-h,x)) - g(u(t+t_0-h,x))\leq
0, $$
for all $
 (t,x)\in[0,h]\times \R\setminus\{L_-\cup L_+\}, 
$
the maximum principle assures that  the function $\delta(t,x)$
 in $R[r,h]$ is ether negative or it reaches a non-negative maximum 
at a point $P_1= (t_1,x_1)$ belonging to  $\partial R_{1}\cup\partial R_{+}\cup\partial
R_{-}\backslash\{h\}\times(-r,r)$. It is easy to see that $P_1 \not\in L_\pm$. Indeed, if 
$P_1 \in L_\pm$ (see Fig. 1) then $\delta_x(P_1+)- \delta_x(P_1-) = \gamma \lambda e^{-\gamma t_1}>0$. Thus the non-negative maximum of $\delta(t,x)$ on $R[r,h]$  is
attained at a point from the parabolic boundary of $R[r,h]$. 
In consequence,   the usual maximum principle holds for each $R[r,h]$ so that,  just as it was done in  Step V of the proof of Lemma  \ref{mlem},  we can appeal to the 
Phragm\`en-Lindel\"of principle in order to conclude that  $\delta(t,z)
\leq 0$ for  all $t \in [0,h], \ z \in \R$.  Applying the above argument consecutively on the intervals $[h,2h],$ $[2h, 3h], \dots$ we find that $\delta(t,x) \leq 0$ for all $t \geq -h, x \in \R$.  Therefore, in view of (\ref{eq14.2}), 
$$
{u}(t+t_0,0) \geq 2{\phi}( {c}t-\epsilon(\infty))-\kappa- \gamma e^{-\gamma
t}\geq \kappa-q'e^{-{\gamma}t},\quad   t \geq -h, 
$$
for some sufficiently large  $q' >\gamma$. Obviously, this yields (\ref{eq14}) with  appropriate $q > q'$.  \hfill \ter
\begin{corollary}\label{x0c} The conclusion of Lemma \ref{x0}  holds without the assumption  ${\lambda_1}<-{\lambda}_{3}$. \end{corollary}
{\it Proof.}  First, we observe that there exists  a monotone function $\hat g(x) \leq g(x)$ satisfying the hypothesis {\rm \bf(H)}  and such that the equation 
\begin{eqnarray} \label{eh1}
u_{t}(t,x) = u_{xx}(t,x) - u(t,x) + \hat g(u(t-h,x))
\end{eqnarray}
has a pushed  wavefront $\hat \phi(\hat ct+x)$ with the 
associated eigenvalues $\hat \lambda_1$, $\hat \lambda_3= \lambda_3$ such that  $\hat \lambda_1 < - \hat \lambda_3$.  Indeed, let  $g_n(x)\leq g(x)$,  be a sequence of monotone functions satisfying {\rm \bf(H)}, coinciding with $g(x)$ on $[1/n, \kappa]$, uniformly on $[0,\kappa]$ converging to $g(x)$ and such that $\lim_{n\to +\infty} g_n'(0) =1$. Then \cite[Lemma 3.5]{LZh} implies that $c_n:= c_*(g_n) \leq  c := c_*(g)$ while the proof of Proposition \ref{mainP} shows that 
$\liminf_{n \to +\infty} c_n \geq c$. This means that $\lim_{n \to +\infty} c_n = c > c_\#> c^{(n)}_\#$ and $\lim_{n \to +\infty} \lambda^{(n)}_1 =0 <  -\lambda_3$ where, similarly  to $c_\#, \lambda_1$, the numbers  $c^{(n)}_\#, \lambda^{(n)}_1$ are determined  from the characteristic equation (\ref{cra}) with $g'(0)$  replaced by $g'_n(0)$. 

In consequence, if $\hat u(t,x)$ denotes the solution of 
the initial value problem (\ref{e2}) for (\ref{eh1}), with $w_0\not\equiv 0$, then Lemma \ref{x0} implies that $\hat u(t,0)  \geq \kappa - qe^{-\nu t}$, $t >0,$ for some positive $q, \nu$. Finally, by comparing initial value problems (\ref{e1}), (\ref{e2}) and (\ref{eh1}), (\ref{e2}) and invoking the Phragm\`en-Lindel\"of principle, we get that $u(t,0) \geq \hat u(t,0) \geq \kappa - qe^{-\nu t}$ for all $t >0$.
 \hfill \ter
\begin{corollary} \label{bs}
Assume that all the conditions of Theorem \ref{main2} are
satisfied. Then there exist $K >1, t_1 >0$ and $z',z'' \in \R$ such that
\begin{eqnarray*} \hspace{-25mm}
u(t,x) \geq \phi(-|x|+ct-z')-Ke^{-\gamma t}\eta(-|x|+ct-z'') \ \mbox{for all} \  t >t_1-h,\ x \in \R. \end{eqnarray*}
\end{corollary}
{\it Proof. } Consider
$
u_{-}(t,x)=\phi(-|x|+ct-\epsilon(t))-\gamma e^{-\gamma t}\theta(t,x). 
$
Analysing the proof of Claim I of Lemma \ref{x0}, we can easily find  that it is also valid for $x\not=0$
if we replace ${\phi}_{+}(t,x)$ with $\kappa$. 
Moreover,  in such a case, the restriction 
${\lambda_1}<-{\lambda}_{3}$ is unnecessary (recall that this restriction appears due to 
the term ${\phi}_{+}(t,x) -\kappa$).  Hence, we  conclude that,  for an
appropriate choice of $\epsilon(t)$,  it holds $\mathcal{N}_1u_{-}(t,x) \leq 0$ for all $x > 0, \ t >0,\ (t,x) \not\in L_+$. 
Since $\mathcal{N}_1u_{-}(t,x) = \mathcal{N}_1u_{-}(t,-x)$, we conclude  that  $u_{-}$ is a sub-solution in the 
region $x\not=0, \ t >0,\ (t,x) \not\in L_\pm$.  In addition, for some sufficiently large $t_1>0$, it holds
$$
u(t+t_1,0) \geq \kappa - {q}e^{-{\gamma}(t+t_1)} > \phi(ct-\epsilon(-h))-\gamma e^{-\gamma t}\eta(ct-\epsilon(\infty)-z_{1}) = 
$$ $$\phi(ct-\epsilon(-h))-\gamma e^{-\gamma t}
>  u_{-}(t,0) \ \mbox{for all} \ t \geq -h. 
$$
Now, arguing as in Claim II of the proof of Lemma \ref{x0}, we can also assume that $t_1$ is chosen in such a way that $u_{-}(s,x)\leq u(s+t_1,x)$ for $x \in \R, \ s \in [-h,0]$. But then, using the Phragm\`en-Lindel\"of principle
in the regions $[hj,h(j+1)] \times [0,+\infty)$, $[hj,h(j+1)] \times (-\infty,0] $, $ j = 0,1, \dots $ according to the procedure established in the last paragraph of  the proof of Lemma \ref{x0}, we conclude that, for all $x \in \R, \ t \geq t_1-h$, it holds that 
$$
u(t,x)  \geq u_{-}(t-t_1,x) = \phi(-|x|+c(t-t_1)-\epsilon(t-t_1))-\gamma e^{-\gamma (t-t_1)t}\theta(t-t_1,x)\geq $$
$$
\phi(-|x|+c(t-t_1)-\epsilon(\infty))-\gamma e^{-\gamma (t-t_1)t}\eta(-|x|+{c}(t-t_1)-\epsilon(\infty)-z_{1}).
$$
This completes the proof of the corollary. \hfill \ter

\begin{lem}\label{lem4}
 Assume all the conditions of  Theorem \ref{main2} and suppose that for some sequence $t_n \to +\infty$ and $s_1,s_2 \in \R$, it holds 
 \begin{eqnarray} 
 \label{541}
\lim_{n \to \infty }\sup_{x\leq 0} |u(t_n+s,x) -
\phi(x +c(t_n+s)+ s_1)|/\eta(x+ct_n)=0,\\ \label{542}
\lim_{n \to \infty }\sup_{x\geq
0}|u(t_n+s,x) -   \phi(-x +c(t_n+s)+ s_2)|/\eta(-x+ct_n) =0,
\end{eqnarray}
uniformly on 
$s\in [-h,0]$.  Then for every $\delta>0$ there exists  $T(\delta)>0$ such that
\begin{eqnarray}\label{stab}
\sup_{x\leq 0}\frac{|u(t,x)-\phi(x+ct+s_{1})|}{\eta(x+ct)}<\delta \quad \mbox{for all} \  t\geq T(\delta), 
\end{eqnarray}
$$
\sup_{x\geq 0}\frac{|u(t,x)-\phi(-x+ct+s_{2})|}{\eta(-x+ct)}<\delta \quad \mbox{for all} \  t\geq T(\delta). 
$$
\end{lem}
{\it Proof.}  It suffices to establish (\ref{stab}), since $u(t,-x)$ also solves equation (\ref{e1}) and satisfies  all the hypotheses of Theorem \ref{main2}.  Without restricting generality, we  can take $s_1=0$. 
We know from Corollary \ref{x0c} that $
u(t,0)\geq \kappa-qe^{-\nu t}, \  t\geq 0. 
$
Fix $\gamma\in (0,\min\{\nu, -c\lambda_3\})$  and consider $\epsilon(t)=\alpha\delta\gamma^{-1}e^{-\gamma t}$ (with $\alpha$ defined in Step II of  Lemma \ref{lem1}) and 
\begin{eqnarray*}
{u}_n(t,x) = \phi(x+ct+ct_n- \alpha\gamma^{-1}e^{\gamma h}\delta+\epsilon(t))-\delta e^{-\gamma
t}\eta(x+ct + ct_n). 
\end{eqnarray*}
Let  positive integer $N=N(\delta)$ be such  that  $\delta e^{\nu t_{N}} >q$ and
\begin{eqnarray*}\hspace{-2cm}
\sup_{(s,x)\in[-h,0]\times(-\infty,0]}\frac{|u(t_{n}+s,x)-\phi(x+c(t_{n}+s))|}{\eta(x+c(t_{n}+s))}<\delta
\ \mbox{for all} \ n\geq N(\delta). 
\end{eqnarray*}
Then we obtain, for all for $(s,x)\in [-h,0]\times(-\infty,0],$ 
\begin{eqnarray}\label{inq02}\nonumber  \hspace{-5mm}
{u}_N(s,x)\leq \phi(x+c(t_{N}+s))-\delta \eta(x+c(t_{N}+s)) \leq
u(t_{N}+s,x).
\end{eqnarray}
Let us show now that a similar relation holds for all $(t,x)\in [t_N,\infty)\times\{0\}$ once $N(\delta)$ is large. 
Indeed, we have that  
$
{u}_N(t,0)\leq\kappa- \delta e^{-\gamma t}$ for all $ t\geq 0
$
so that,   
$$
u(t+t_N,0)-{u}_N(t,0)\geq  \delta  e^{-\gamma t}-qe^{-\nu t_N}
e^{-\nu t}>0, \ t \geq 0. 
$$
Next, observe that ${u}_n(t,x) = w_-(t,x+c(t+t_n))$ where $w_-$ is defined in Lemma \ref{lem1} (by Remark \ref{R22},  the summand $- \alpha\gamma^{-1}e^{\gamma h}$ within the argument of $\phi$ doesn't matter). Since $\delta < \sigma$, we find that   
${(u_n)}_t(t,x)- {(u_n)}_{xx}(t,x) + {u_n}(t,x) - g({u_n}(t-h,x)) = 
({\mathcal N}w_-)(t,x+c(t+t_n))  < 0$  for all 
$(t,x)\in[0,\infty)\times\R$, $x+ct + ct_n\not=0$.  Furthermore, if $x'+ct' +ct_n=0$ at some point $(x',t')$ then ${(u_n)}_x(t',x'+0)- {(u_n)}_x(t',x'-0) = \lambda \delta e^{-\gamma t'} >0$.  Therefore, repeatedly applying  the Phragm\`en-Lindel\"of principle
in the regions $[hj,h(j+1)] \times (-\infty,0] $, $ j = 0,1, \dots $ according to the procedure established in the last paragraph of  the proof of Lemma \ref{x0}, we conclude that, for all $x \leq 0, \ t \geq -h$, 
$$
u(t+t_N,x)\geq{u}_N(t,x)\geq
\phi(x+c(t+t_N)- \alpha\gamma^{-1}e^{\gamma h}\delta)-\delta \eta(x+c(t+t_N)). 
$$ 
Hence, taking positive constant $K_1= K_1(\alpha, \gamma, h)$  as in  (\ref{it}), we obtain that 
\begin{equation}\label{pr1}
\hspace{-5mm}u(t,x)\geq
\phi(x+ct)-\delta (1+K_1) \eta(x+ct), \ t \geq t_N-h, \ x \leq 0.  
\end{equation}
On the other hand, by our assumptions,  for all for $(s,x)\in [-h,0]\times(-\infty,0],$  
\begin{equation}\label{te}
{u}(t_{N}+s,x)\leq \phi(x+c(t_{N}+s))+\delta\eta(x+c(t_{N}+s)).
\end{equation}
If, in addition, $N=N(\delta)$ is so large that 
$$
\phi(c(t_{N}+s))+\delta\eta(c(t_{N}+s)) >\kappa,\ s \in [-h,0],$$ then  (\ref{te})  holds also for all $(s,x)\in [-h,0]\times \R$.   
Therefore for $\delta \in (0,q_0]$, by  Lemma \ref{lem1}, 
$$
u(t+t_N,x) \leq \phi(x+c(t_N+t)+C\delta)+\delta e^{-\gamma
t}\eta(x+c(t_N+t)), \ t \geq 0, \ x \in \R, 
$$
for positive $C>0$ defined in  Lemma \ref{lem1}.   Next, due to (\ref{it}),  for all $(t,x) \in \R^2$, we have 
$$
\phi(x+c(t_N+t)+  C\delta) \leq
\phi(x+c(t_N+t))  + K_1\delta \eta(x+c(t_N+t)).
$$
In consequence, we obtain 
$$
u(t,x)\leq\phi(x+ct)+\delta(1+K_1)\eta(x+ct), \ \mbox{for} \ t \geq t_N, \ x \leq 0. 
$$
The latter inequality together with  (\ref{pr1}) imply (\ref{stab}). 
 \hfill \ter

\section{Proof of Theorem \ref{main2}: main arguments}\label{sub4}

Set $z=x+ct$ and $w(t,z):= u(t,x)=u(t,z-ct)$, then $w(t,z)$
satisfies equation (\ref{wz}), (\ref{wzic}) for $(t,z) \in \R_+\times \R$ and
possesses a compact and invariant $\omega$-limit set $
\omega(w_{0})$ defined in Remark \ref{ols}. Consider the semi-infinite strip $ \Omega=\{(s,z)\in[-h,0]\times\R,\ z\leq -ch\}. $
By Corollary \ref{coro1} and Remark \ref{dopi}, for some $K>0,\
\zeta_1 \in \R,$  it holds
\begin{eqnarray}\label{edo1}\hspace{0mm}
w(t,z)\leq \phi(z+\zeta_{1})+Ke^{-\gamma t}\eta(z+\zeta_{1}),
\quad   z\in \R, \  t\geq-h.
\end{eqnarray}
Therefore  the  set
$$
A = \{a\in \R: v(s,z)\leq \phi(z+a), \ (s,z) \in \Omega, \
\mbox{for each } \ v \in  \omega(w_0)\}
$$
in non-empty. Since, by Corollary \ref{bs},
\begin{eqnarray}\label{edo}
\phi(z-z_{1})-K\gamma e^{-\gamma t}\eta(z-z_{1})\leq w(t,z), \quad
z\leq ct,\   t\geq-h,
\end{eqnarray}
$A$ is bounded below. Set $\hat{a}:=\inf A$, obviously,
$\hat{a}\in A.$  We claim that $v_*(s_*,z_*)=\phi(z_*+\hat{a})$
for some $(s_*,z_*)\in\Omega$ and  $ v_{*}\in\omega(w_{0})$.
Indeed, suppose on the contrary that
\begin{eqnarray}\label{inq1}
  v(s,z)<\phi(z+\hat{a})\ \mbox{for all} \ (s,z)\in\Omega, \
 v\in\omega(w_{0}).
\end{eqnarray}
For positive $\varsigma$ and an entire 
solution  $v\in\omega(w_0), \ v: \R^2 \to [0,\kappa]$, consider
$\rho(t,z)=v(t,z)-\phi(z+\hat{a}-\varsigma)$. 
Let  $R>ch$ be such that $\phi(-R+\zeta_{1})<\delta$.  Then, for each  $\xi(t,z)$ lying 
between points  $v(t-h,z-ch)$ and $\phi
(z+\hat{a}-\varsigma-ch)$ with $z\leq -R, \ t \in \R$, we have  $\xi(t,z) \in (0, \delta)$.  
Next, set  $r(t,z)=\eta(z)e^{-\gamma t}$ and let $\delta$ be as in
(\ref{C1a}).  In view of (\ref{C1a}),
\begin{eqnarray*}\label{eq3.1} && \hspace{-25mm}
r_{t}(t,z)-r_{zz}(t,z)+cr_{z}(t,z)+r(t,z)= \eta(z)e^{-\gamma
t}[1-\gamma-\lambda^{2}+c\lambda]\\\nonumber && \hspace{0mm} \geq
\eta(z)e^{-\gamma t}g'(\xi(t,z))e^{-\lambda ch+\gamma h} ,\quad t
>0, \ z \leq -R, \ \varsigma >0,\ v\in\omega(w_{0}), 
\end{eqnarray*}
\begin{eqnarray*} \hspace{-25mm}
\rho_{t}(t,z)=\rho_{zz}(t,z)-c\rho_{z}(t,z)-\rho(t,z)+g'(\xi(t,z))\rho(t-h,z-ch),
\quad t \in \R,\ z \leq -R.
\end{eqnarray*}
On the other hand,  since the set $\omega(w_{0})$ is  compact and
invariant (the latter means that $\omega(w_{0})$  consists of {\it entire} solutions $v: \R^2 \to [0,\kappa]$) and $\phi$ increases on $\R$,  we can fix $\varsigma>0$
such that $(\ref{inq1})$ implies
\begin{equation}\label{eq3.2} \hspace{-25mm}
v(t,z)< \phi(z+\hat a-\varsigma),\ t>0,\ -R\leq z\leq -ch, \ v \in
\omega(w_0).
\end{equation}
Without loss of generality, we can also suppose that $\varsigma$ is
sufficiently small to meet
\begin{eqnarray}\label{eq3.3} \hspace{-15mm}
\phi(z+\hat{a})<\phi(z+\hat{a}-\varsigma)+\eta(z)
e^{-\gamma s}\quad  \mbox{for all}\ z\in\R,\ s \in[-h,0].
\end{eqnarray}
Now, we set $\delta(t,z):= \rho(t,z) - r(t,z)$. Note that, by
(\ref{inq1}) and (\ref{eq3.3}),  for all for $s\in [-h,0]$, $z\leq -ch$, it holds
 $$\delta(s,z)= v(s,z)-(\phi(z+\hat{a}-\varsigma) +\eta(z)e^{-\gamma s}) < v(s,z)-\phi(z+\hat{a}) < 0,$$
and therefore, in virtue of  the above mentioned  properties of $\rho, r$,
\begin{eqnarray*} && \hspace{-25mm}
\delta_{zz}(t,z)- \delta_{t}(t,z)-c\delta_{z}(t,z)-\delta
(t,z)\geq -g'(\xi(t,z))\rho(t-h,z-ch)\\\nonumber &&\hspace{-25mm}
+\eta(z)e^{-\gamma t}g'(\xi(t,z))e^{-\lambda ch+\gamma h}
=-g'(\xi(t,z))\delta(t-h,z-ch)> 0 \  \mbox{for}\ z\leq -R,\ t
\in [0,h].
\end{eqnarray*}
Taking into account that,  due to  (\ref{eq3.2}), it holds
$-\kappa -1 < \delta(t,z) < 0$  for all $t\in [0,h],\ -R \leq z \leq -ch,$ we
can invoke now the Phragm\`en-Lindel\"of principle \cite{PW} in
order to conclude that $\delta(t,z) < 0$ for all  $t\in [0,h],\
z \leq -R$. But then, by repeating the above argument for the time
intervals $[h,2h], \ [2h, 3h], \dots$, and using (\ref{eq3.2}) we
conclude that
$$v(t,z)\leq \phi(z+\hat{a}-\varsigma) +\eta(z)e^{-\gamma t}$$
for all $t \geq 0, \ z \leq -ch$.  Due to the invariance property  of  $\omega(w_0)$ this yields
$$v(s,z)< \phi(z+\hat a-\varsigma),\ -h\leq s\leq 0,\  z \leq -ch, \ v\in  \omega(w_0),$$
contradicting the definition of $\hat a$.

Hence, $w_{*}(s_*,z_*)=\phi(z_*+\hat{a})$ for some
$(s_*,z_*)\in\Omega$ and   $w_{*}\in\omega(w_{0})$. Therefore,  by the
strong principle maximum and invariance property of
$\omega(w_{0})$, we obtain  that
 $\phi \in\omega(w_{0}).$

Next, it follows from  (\ref{edo1}) and (\ref{edo}) that, for all  $z\leq ct, \  t\geq-h$, it holds 
$$
 |w(t,z)- \phi(z+\hat a)| \leq 
 \phi(z+\zeta_{1})-  \phi(z-z_{1})+Ke^{-\gamma t}(\eta(z-z_{1}) +\eta(z+\zeta_{1})). 
$$ 
In consequence, for each  $\epsilon>0$ we can find 
$T(\epsilon)>0$ such that
$$
|w(t+s,z)-\phi(z+\hat{a})|<\epsilon\quad\mbox{for}\ t \geq T(\epsilon),\  cT(\epsilon) \leq   z \leq ct, \ s \in [-h,0].  
$$
 and
$$
\frac{|w(t,z)-\phi(z+\hat{a})|}{\eta(z)}<\epsilon\quad\mbox{for} \
 t \geq T(\epsilon),\   z \leq -cT(\epsilon), \ s \in [-h,0]. 
$$
On the other hand, since  $\phi\in\omega(w_0)$, there exist $t_n\rightarrow\infty$ and an integer 
$n(\epsilon)$ so that:
$$
\frac{|w(t_n+s,z)-\phi(z+\hat{a})|}{\eta(-cM)}<\epsilon, \quad n\geq
n(\epsilon),\ |z|\leq cT(\epsilon), \ s \in [-h,0]. 
$$
Obviously, the last three inequalities imply (\ref{541}).
 Moreover,  by considering the solution $\hat{u}(t,x) = u(t,-x)$ together with  the obtained sequence $\{t_n\}$,
  we can see that (\ref{542}) is also satisfied for a subsequence $\{t_{n_j}\}\subset \{t_n\}$ and an appropriate $s_2$.
 Finally, an application of  Lemma \ref{lem4} completes the proof of Theorem  \ref{main2}.

\section*{Acknowledgments}  
This research was supported by FONDECYT (Chile).

\end{document}